# A CONJUGATE PRIOR FOR DISCRETE HIERARCHICAL LOG-LINEAR MODELS

By Hélène Massam,[1] Jinnan Liu and Adrian Dobra

*York University, York University and University of Washington*

In Bayesian analysis of multi-way contingency tables, the selection of a prior distribution for either the log-linear parameters or the cell probabilities parameters is a major challenge. In this paper, we define a flexible family of conjugate priors for the wide class of discrete hierarchical log-linear models, which includes the class of graphical models. These priors are defined as the Diaconis–Ylvisaker conjugate priors on the log-linear parameters subject to "baseline constraints" under multinomial sampling. We also derive the induced prior on the cell probabilities and show that the induced prior is a generalization of the hyper Dirichlet prior. We show that this prior has several desirable properties and illustrate its usefulness by identifying the most probable decomposable, graphical and hierarchical log-linear models for a six-way contingency table.

**1. Introduction.** We consider data given under the form of a contingency table representing the classification of $N$ individuals according to a finite set of criteria. We assume that the cell counts in the contingency table follow a multinomial distribution. We also assume that the cell probabilities are modeled according to a hierarchical log-linear model. The class of discrete graphical models that are Markov with respect to an arbitrary undirected graph $G$ is an important subclass of the class of hierarchical log-linear models, in part because graphical models can be interpreted in terms of conditional independences that can easily be read off of the graph, and in part because they allow for parsimony in the number of parameters in the analysis of complex high-dimensional data. We will therefore give special attention to the class of graphical models throughout the paper.

Received March 2008; revised November 2008.

[1]Supported by NSERC Discovery Grant A8946.

In the Bayesian analysis of contingency tables, the selection of a prior distribution for either the log-linear parameters or the cell probabilities parameter is a major challenge (see Clyde and George [3]). For decomposable graphical models, Dawid and Lauritzen [8] have identified a standard conjugate prior which they called the hyper Dirichlet. The hyper Dirichlet presents the mathematical convenience of a conjugate prior; it has the flexibility given by a number of hyperparameters (as many as there are free cell probabilities in the model) and, additionally, has the strong hyper Markov property. The latter is very desirable, since it allows for local updates within prime components, thus simplifying the computation of Bayes factors in a model selection process. For decomposable models, with the hyper Dirichlet as a prior, Bayes factors can be computed explicitly. The hyper Dirichlet has therefore been used in many studies (see, e.g., Madigan and Raftery [21] or Madigan and York [23]). However, it has the disadvantage of being defined only for the class of decomposable graphical models, which, as the number of factors increases, becomes a smaller and smaller part of the class of graphical models and, even more so, of hierarchical models.

Considerable efforts have been devoted to the study of alternative priors valid for the larger class of hierarchical models. Knuiman and Speed [18], Dellaportas and Forster [9] and King and Brooks [17] propose various versions of a multivariate normal prior for the log-linear parameters.

In this paper, we propose a new prior for the class of hierarchical log-linear models. This prior is the Diaconis–Ylvisaker conjugate prior for the log-linear parameters subject to baseline constraints. We show that it is a generalization of the hyper Dirichlet to nondecomposable graphical models and, even more generally, to hierarchical log-linear models. We also show that, like the hyper Dirichlet, it has the advantage of being a conjugate prior while offering flexibility through its hyperparameters. We illustrate its applicability for the well-known Czech Autoworkers example previously analyzed by Edwards and Havranek [14], Madigan and Raftery [21] and Dellaportas and Forster [9]. We employed MC$^3$ to explore the space of decomposable, graphical and hierarchical log-linear models for this six-dimensional binary table. Dobra and Massam [13] and Dobra et al. [11] demonstrate that our conjugate priors scale to higher-dimensional examples arising from social studies as well as gene expression and genomewide association studies.

A secondary aim of this paper is to contribute to a discussion on a question asked by Gutiérrez-Pena and Smith [15], itself motivated by a characterization given, for univariate natural exponential families (henceforth abbreviated NEF) by Consonni and Veronese [5]. The latter proved that the prior induced onto the mean parameter of an NEF, from the Diaconis–Ylvisaker conjugate prior for the log-linear parameters, is standard conjugate if and only if the variance function of the NEF is quadratic in the mean. Leucari [20] showed that, in the case of a decomposable graphical model, the induced prior on the mean parameter is standard conjugate even though the variance function is not quadratic, thus providing a negative answer to the question posed by Guitterez-Pena and Smith [15] as to whether the characterization of Consonni and Veronese [5] could be extended to multivariate NEF. Here, we show more precisely that the induced prior on the clique and separator marginal probabilities, which we will denote $p^G$, is standard conjugate (it is the hyper Dirichlet as mentioned above) while the induced prior on the cell probabilities parametrization, denoted $p_{\mathcal{D}}$, does not share the same property. This is achieved through the derivation of the prior induced on the cell probabilities from our prior on the log-linear parameters.

The paper is organized as follows. In Section 2, we define the parameters we chose to use in order to express the multinomial distribution. They are the classical log-linear parameters defined by the "baseline" or "corner" constraints, and we show that this is the parametrization obtained if we make a change of variable from the cell counts to the marginal cell counts. In Section 3, we derive the Diaconis and Ylvisaker [10] (henceforth abbreviated DY) conjugate prior for this parametrization for hierarchical log-linear models, and we study its properties. We first characterize the set of hyperparameters for which our conjugate prior is proper. We then use this characterization to construct a set of hyperparameters that leads to a proper prior. We compute the moments of the prior cell probabilities that can be used to guide our choice of hyperparameters. Finally, we show that, like the hyper Dirichlet, our prior on the log-linear parameters has what we might call the strong hyper Markov property extended to nondecomposable models, so that inference in a Bayesian framework can be made prime component by prime component. In Section 4, we derive the induced prior for the cell probabilities for the decomposable graphical model, the nondecomposable graphical model and, more generally, the general log-linear hierarchical model. As mentioned above, we discuss the conjecture of Gutiérrez-Pena and Smith [15]. In Section 5, we present a comprehensive analysis of the Czech Autoworkers data that includes a sensitivity study about the influence of the conjugate prior specification on the highest posterior probability log-linear models. In Section 6 we briefly talk about additional developments using our conjugate prior. Major proofs of some of our results herein are given in the Appendix.

## 2. The log-linear model.

2.1. *The parametrization.* Let $V$ be the set of criteria. Let $X = (X_\gamma, |\gamma \in V)$ such that $X_\gamma$ takes its values (or levels) in the finite set $I_\gamma$ of dimension $|I_\gamma|$. When a fixed number of individuals are classified according to the $|V|$ criteria, the data is collected in a contingency table with cells indexed

by combination of levels for the $|V|$ variables. We adopt the notation of Lauritzen [19] and denote a cell by

$$i = (i_\gamma, \gamma \in V) \in \mathcal{I} = \underset{\gamma \in V}{\times} \mathcal{I}_\gamma.$$

The count in cell $i$ is denoted $n(i)$, and the probability of an individual falling in cell $i$ is denoted $p(i)$. For $E \subset V$, cells in the $E$-marginal table are denoted $i_E \in \mathcal{I}_E = \underset{\gamma \in E}{\times} \mathcal{I}_\gamma$ and the marginal counts are written

$$(2.1) \qquad n(i_E) = \sum_{j \in \mathcal{I}_{V \setminus E}} n(i_E, j_{V \setminus E}).$$

For $N = \sum_{i \in \mathcal{I}} n(i), (n) = (n(i), i \in \mathcal{I})$ follows a multinomial $\mathcal{M}(N, p(i), i \in \mathcal{I})$ distribution with probability density function

$$(2.2) \qquad P((n)) = \binom{N}{(n)} \prod_{i \in \mathcal{I}} p(i)^{n(i)}.$$

Let $i^*$ be a fixed but arbitrary cell which, for convenience, we take to be the cell indexed for each factor by the "lowest level" itself indexed, for convenience again, by 0. Thus, $i^*$ is the cell

$$i^* = (0, 0, \dots, 0).$$

We now have to choose a parametrization for the log-linear model; that is, a parametrization for $\log p(i)$. As shown in Darroch and Speed [7] (see also Lauritzen [19], Appendix B.2), each possible metric in the space $\mathbb{R}^{\mathcal{I}}$ of real-valued functions defined on $\mathcal{I}$ corresponds to a different parametrization of the log-linear model. Moreover, as illustrated in Wermuth and Cox [29], a given parametrization can be best suited to a given type of problem. There is, therefore, no "best" parametrization in general.

In this paper, we choose to work with the parametrization given by "baseline" or "corner" constraints; that is, the parametrization that follows if we choose the "substitution weight" metric for the space $\mathbb{R}^{\mathcal{I}}$, as given in Section 3.1 of Darroch and Speed [7]. This parametrization has the practical advantage of yielding the marginal counts as the canonical statistic in the exponential family form of (2.2), thus making the derivation of the general form of marginal and conditional distributions, as well as that of conjugate distributions, very easy to express. The log-linear parameters are

$$(2.3) \qquad \theta_E(i_E) = \sum_{F \subseteq E} (-1)^{|E \setminus F|} \log p(i_F, i^*_{F^c}),$$

which, by Moebius inversion, is equivalent to

$$(2.4) \qquad p(i_E, i^*_{E^c}) = \exp \sum_{F \subseteq E} \theta_F(i_F).$$

We note that $\theta_\varnothing(i^*) = \log p(i^*), i \in \mathcal{I}$ and we will therefore adopt the notation

$$(2.5) \qquad \theta_\varnothing(i^*) = \theta_\varnothing, \qquad p(i^*) = p_\varnothing = \exp \theta_\varnothing.$$

The parametrization (2.3) was first used by Mantel [24]. It is used in most standard statistical software such as GLIM or R (see Agresti [1], page 150). It has recently been used in Consonni and Leucari [4], in the case of binary data. It seems, however, that it is less commonly used in the literature than the so-called $u$-parametrization (see Bishop, Fienberg, and Holland [2]) though the "interaction" terms $\theta_E(i_E)$ in (2.3) are easy to interpret as ratios of log-odds ratios or as partial cross-product ratios. Indeed, one can easily verify (see Lauritzen [19], page 37) that for any $\alpha, \beta$ in $E$, with the notation $E^- = E \setminus \{\alpha, \beta\}$, $\theta_E(i_E)$ can also be written as the alternating sum of conditional log-odds ratios

$$\theta_E(i_E) = \sum_{F \subseteq E^-} (-1)^{|E \setminus F|} \log \frac{p(i_\alpha, i_\beta | i_F, i^*_{(E^- \setminus F)}) p(i^*_\alpha, i^*_\beta | i_F, i^*_{(E^- \setminus F)})}{p(i^*_\alpha, i_\beta | i_F, i^*_{(E^- \setminus F)}) p(i_\alpha, i^*_\beta | i_F, i^*_{(E^- \setminus F)})}.$$

Another pleasant feature of this parametrization is that, as we shall see more precisely at the beginning of the next subsection, the parameters given in (2.3) are obtained as the canonical parameters of the multinomial distribution when we make the change of variable from the cell counts to the marginal cell counts. In order to identify which ones of the $\theta_E(i_E)$ defined in (2.3) is a free parameter, we need the following lemma.

LEMMA 2.1.    *If for $\gamma \in E, E \subseteq V$ we have $i_\gamma = i^*_\gamma = 0$, then $\theta_E(i_E) = 0$.*

PROOF.    By definition, and since $(i_{F \cup \gamma}, i^*_{(F \cup \gamma)^c}) = (i_F, i^*_{F^c})$ if $i_\gamma = i^*_\gamma = 0$, we have

$$\theta_E(i_E) = \sum_{F \subseteq E \setminus \gamma} (-1)^{|E \setminus F|} \log p(i_F, i^*_{F^c})$$

$$- \sum_{F \subseteq E \setminus \gamma} (-1)^{|E \setminus F|} \log p(i_{F \cup \gamma}, i^*_{(F \cup \gamma)^c})$$

$$= \sum_{F \subseteq E \setminus \gamma} (-1)^{|E \setminus F|} \log p(i_F, i^*_{F^c})$$

$$- \sum_{F \subseteq E \setminus \gamma} (-1)^{|E \setminus F|} \log p(i_F, i^*_{F^c})$$

$$= 0. \qquad\qquad \square$$

From this lemma, it follows immediately that our parametrization is indeed the "baseline" or "corner" constraint parametrization that sets to 0 the values of the $E$-interaction log-linear parameters when at least one index in $E$ is at level 0 (see Agresti [1], page 150). Therefore, for each $E \subseteq V$, there are only $\prod_{\gamma \in E}(|\mathcal{I}_\gamma| - 1)$ parameters and for any $E \subseteq V$, we introduce the convenient notation

$$(2.6) \qquad \mathcal{I}_E^* = \{i_E | i_\gamma \neq i_\gamma^*, \forall \gamma \in E\}.$$

In words, $\mathcal{I}_E^*$ is the set of marginal cells $i_E$ such that none of their components is equal to 0. We set $\mathcal{I}_V^* = \mathcal{I} \setminus \{i^*\}$. For example, if $E = \{a, b, c\}$, $a$ takes the values $\{0, 1, 2, 3\}$, $b$ takes the values $\{0, 1, 2\}$, $c$ takes the values $\{0, 1\}$, then

$$\mathcal{I}_E^* = \{(1, 1, 1), (2, 1, 1), (3, 1, 1), (1, 2, 1), (2, 2, 1), (3, 2, 1)\}.$$

It will also be convenient to introduce the notation

$$(2.7) \qquad \mathcal{E}_\ominus = \{E \subseteq V, E \neq \varnothing\}$$

for the power set of $V$ deprived of the empty set and the notation $\mathcal{E}$ for the power set of $V$.

By (2.4) and (2.5), we have

$$
\begin{aligned}
(2.8) \quad p_\varnothing &= 1 - \sum_{i \in \mathcal{I}, i \neq i^*} p(i) \\
&= 1 - \sum_{E \in \mathcal{E}_\ominus} \sum_{i_E \in \mathcal{I}_E^*} p(i_E, i_{E^c}^*) = 1 - \sum_{E \in \mathcal{E}_\ominus} \sum_{i_E \in \mathcal{I}_E^*} \exp \sum_{F \subseteq E} \theta_F(i_F) \\
&= 1 - \sum_{E \in \mathcal{E}_\ominus} \sum_{i_E \in \mathcal{I}_E^*} \exp\left(\theta_\varnothing + \sum_{F \subseteq E, F \neq \varnothing} \theta_F(i_F)\right) \\
&= 1 - p_\varnothing \sum_{E \in \mathcal{E}_\ominus} \sum_{i_E \in \mathcal{I}_E^*} \exp\left(\sum_{F \subseteq E, F \neq \varnothing} \theta_F(i_F)\right).
\end{aligned}
$$

In order to simplify our notation, from now on, we will use

$$F \subseteq_\ominus E$$

to express that $F$ is included in E but is not equal to the empty set and, for $i_E \in \mathcal{I}_E^*$, $E \in \mathcal{E}$, the notation

$$i(E) = (i_E, i_{E^c}^*),$$

which is not to be confused with $i_E$, the $E$ marginal cell. We will also write $\theta(i_E)$ for $\theta_E(i_E)$. Then, (2.8) yields

$$(2.9) \quad p(i(E)) = \frac{\exp \sum_{F \subseteq_\ominus E} \theta(i_F)}{1 + \sum_{E \in \mathcal{E}_\ominus} \sum_{j_E \in \mathcal{I}_E^*} \exp(\sum_{F \subseteq_\ominus E} \theta(j_F))}, \qquad E \in \mathcal{E}.$$

We note, in particular, that

$$(2.10) \qquad p_\varnothing = \frac{1}{1 + \sum_{E \in \mathcal{E}_\ominus} \sum_{j_E \in \mathcal{I}_E^*} \exp(\sum_{F \subseteq_\ominus E} \theta(j_F))}.$$

2.2. *The multinomial distribution for discrete data.* We now want to give the probability density function of the multinomial distribution under the form of an exponential family when the statistical model is a hierarchical log-linear model. Let us first show that the parameters in (2.3) are the canonical parameters of the multinomial distribution for the saturated model after we make the change of variables

$$(2.11) \quad (n) = (n(i), i \in \mathcal{I}^*) \mapsto Y = (y(i_E) = n(i_E), E \in \mathcal{E}_\ominus, i_E \in \mathcal{I}_E^*)$$

from joint cell counts to marginal cell counts as defined in (2.1).

Lemma 2.2. *The probability function of the multinomial distribution as given in (2.2) can be represented as a natural exponential family, with canonical parameters* $\theta(i_E), E \in \mathcal{E}_\ominus, i_E \in \mathcal{I}_E^*$ *as defined in (2.3) and with canonical statistics the marginal cell counts* $(n(i_E), E \in \mathcal{E}_\ominus, i_E \in \mathcal{I}_E^*)$, *as follows:*

$$
\begin{aligned}
(2.12) \quad \prod_{i \in \mathcal{I}} p(i)^{n(i)} = \exp\Big\{ &\sum_{E \in \mathcal{E}_\ominus} \sum_{i_E \in \mathcal{I}_E^*} n(i_E)\theta(i_E) \\
&- N\log\Big(1 + \sum_{E \in \mathcal{E}_\ominus} \sum_{i_E \in \mathcal{I}_E^*} \exp\sum_{F \subseteq \ominus E} \theta(i_F)\Big)\Big\}.
\end{aligned}
$$

Proof. We have

$$
\begin{aligned}
\prod_{i \in \mathcal{I}} p(i)^{n(i)} &= p_\varnothing^{n(i^*)} \prod_{E \in \mathcal{E}_\ominus} \prod_{i_E \in \mathcal{I}_E^*} p(i(E))^{n(i(E))} \\
&= p_\varnothing^{n(i^*)} \prod_{E \in \mathcal{E}_\ominus} \prod_{i_E \in \mathcal{I}_E^*} \Big( \exp\sum_{F \subseteq E} \theta(i_F) \Big)^{n(i(E))} \\
&= \prod_{E \in \mathcal{E}_\ominus} \prod_{i_E \in \mathcal{I}_E^*} \exp\Big( n(i(E)) \sum_{F \subseteq \ominus E} \theta(i_F) \Big) p_\varnothing^{n(i^*) + \sum_{E \in \mathcal{E}_\ominus} \sum_{i_E \in \mathcal{I}_E^*} n(i(E))} \\
&= p_\varnothing^N \exp \sum_{E \in \mathcal{E}_\ominus} \sum_{i_E \in \mathcal{I}_E^*} \Big( n(i(E)) \sum_{F \subseteq \ominus E} \theta(i_F) \Big) \\
&= p_\varnothing^N \exp \sum_{E \in \mathcal{E}_\ominus} \sum_{i_E \in \mathcal{I}_E^*} n(i_E)\theta(i_E) \\
&= \exp\Big\{ \sum_{E \in \mathcal{E}_\ominus} \sum_{i_E \in \mathcal{I}_E^*} \theta(i_E) n(i_E) + N\theta_\varnothing \Big\},
\end{aligned}
$$

where the second equality is due to (2.4), the third to (2.5), the fourth to the identification of the exponent of $p_\varnothing$ as the total count $N$, the fifth to (2.1) and the sixth to (2.5) again. Finally, (2.12) follows from (2.10). □

From the change of variable (2.11) and Lemma 2.2, it follows immediately that the family of distributions of $Y$ is the natural exponential family

$$(2.13) \quad \mathcal{F}_\mu = \left\{ f(y;\theta)\mu(y) = \frac{\exp\{\sum_{E \in \mathcal{E}_\ominus} \sum_{i_E \in \mathcal{I}_E^*} \theta(i_E)y(i_E)\}}{(1 + \sum_{E \in \mathcal{E}_\ominus, i_E \in \mathcal{I}_E^*} \exp \sum_{F \subseteq_\ominus E} \theta(i_F))^N} \mu(y), \right.$$
$$\left. \theta \in \mathbb{R}^{\sum_{E \in \mathcal{E}_\ominus} \prod_{\gamma \in E}(|I_\gamma|-1)} \right\},$$

where $\mu$ is a reference measure of no particular interest to us here. This gives us the density for the saturated model.

Let us now consider the hierarchical log-linear model generated by the class $\mathcal{A} = \{A_1, \ldots, A_k\}$ of subsets of $V$, which, without loss of generality, we can assume to be maximal with respect to inclusion. We write

$$(2.14) \quad \mathcal{D} = \{E \subseteq_\ominus A_i \text{ for some } i = 1, \ldots, k\}$$

for the indexing set of all possible interactions in the model, including the main effects. It follows from the theory of log-linear models (see also Darroch and Speed [7]) and from Lemma 2.1 that the model for the cell counts $p(i)$ is the log-linear model with generating class $\mathcal{A}$ if and only if the following constraints are satisfied

$$(2.15) \quad \theta(i_E) = 0, \qquad E \notin \mathcal{D}.$$

Therefore, in this case, for $i_E \in \mathcal{I}_E^*$, (2.4) becomes

$$(2.16) \quad \log \frac{p(i(E))}{p_\varnothing} = \sum_{F \subseteq E, F \in \mathcal{D}} \theta(i_F).$$

Let us now consider an undirected graph $G$ with vertex set $V$. Darroch, Lauritzen and Speed [6] have shown that, for the subclass of graphical models Markov with respect to $G$, the generating class is equal to the set of cliques of $G$, that is, the set of maximal complete subsets of $G$. Therefore, for this subclass,

$$(2.17) \quad \mathcal{D} = \{D \subseteq_\ominus V | D \text{ complete}\}.$$

In general, for the class of hierarchical models with generating class $\mathcal{A}$, the nonzero free log-linear parameters are

$$(2.18) \quad \theta_\mathcal{D} = \{\theta(i_D), D \in \mathcal{D}, i_D \in \mathcal{I}_D^*\}.$$

Let us adopt the short notation

$$F \subseteq_\mathcal{D} D$$

to indicate that $F \subseteq_{\ominus} D$ and $F \in \mathcal{D}$. Then, for the hierarchical log-linear model, (2.9) and (2.10) become

$$(2.19) \quad p(i(E)) = \frac{\exp \sum_{D \subseteq_{\mathcal{D}} E} \theta(i_D)}{1 + \sum_{D \in \mathcal{E}_{\ominus}} \sum_{j_D \in \mathcal{I}_D^*} \exp(\sum_{F \subseteq_{\mathcal{D}} D} \theta(j_F))}, \qquad E \in \mathcal{E}_{\ominus},$$

$$(2.20) \quad p_{\varnothing} = \frac{1}{1 + \sum_{D \in \mathcal{E}_{\ominus}} \sum_{j_D \in \mathcal{I}_D^*} \exp(\sum_{F \subseteq_{\mathcal{D}} D} \theta(j_F))}.$$

Through an argument parallel to that given in Lemma 2.2, it follows that, in the case of a log-linear model with generating class $\mathcal{A}$, the family $\mathcal{F}_{\mu}$ in (2.13) becomes

$$(2.21) \qquad \mathcal{F}_{\mu_{\mathcal{D}}} = \{f_{\mathcal{D}}(y; \theta_{\mathcal{D}}) \mu_{\mathcal{D}}(y), \theta_{\mathcal{D}} \in \mathbb{R}^{d_{\mathcal{D}}}\},$$

where, as in the saturated case, the measure $\mu_{\mathcal{D}}(y)$ is of no particular interest to us here, $\theta_{\mathcal{D}} = (\theta(i_D), D \in \mathcal{D}, i_D \in \mathcal{I}_D^*)$ is the canonical parameter, the dimension $d_{\mathcal{D}}$ of the parameter space is equal to (see Darroch and Speed [7], Proposition 4.3)

$$d_{\mathcal{D}} = \sum_{D \in \mathcal{D}} \prod_{\gamma \in D} (|I_{\gamma}| - 1)$$

and

$$
\begin{aligned}
(2.22) \quad f_{\mathcal{D}}(y; \theta_{\mathcal{D}}) = \exp\Big\{ &\sum_{D \in \mathcal{D}} \sum_{i_D \in \mathcal{I}_D^*} \theta(i_D) y(i_D) \\
&- N \log\Big(1 + \sum_{E \in \mathcal{E}_{\ominus}} \sum_{i_E \in \mathcal{I}_E^*} \exp \sum_{F \subseteq_{\mathcal{D}} E} \theta(i_F)\Big) \Big\}.
\end{aligned}
$$

It is important to note here that, correspondingly to (2.18), only cell probabilities of the form $p(i(D)), D \in \mathcal{D}, i_D \in \mathcal{I}_D^*$, will be free probabilities, since, by Lemma 2.1, all others can be expressed in terms of

$$(2.23) \qquad p_{\mathcal{D}} = (p(i(D)), D \in \mathcal{D}, i_D \in \mathcal{I}_D^*),$$

which will be the cell probability parameter of the multinomial distribution of the hierarchical log-linear model.

When the data is binary, that is, when $|\mathcal{I}_{\gamma}| = 2, \gamma \in V$, there is only one element in $\mathcal{I}_E^*$ for each $E \in \mathcal{E}_{\ominus}$; therefore, since $\theta(i_E)$ is zero if $i_E \notin \mathcal{I}_E^*$, we can use the simplified notation $\theta(E)$, $p(E)$ and $y(E) = n(E)$, respectively, for the canonical parameters $\theta(i_E)$ in (2.3), the cell probabilities $p(i(E))$ in (2.9) and the marginal counts in (2.11) in all the formulas above.

We note that (2.23) becomes $p_{\mathcal{D}} = (p(D), D \in \mathcal{D})$.

Let us illustrate this notation with an example. Let $G$ be the graph with vertices $a, b, c, d$ and edges $(a, b), (b, c), (c, d)$ and $(d, a)$. In this case, the graphical model is actually the same as the hierarchical model with generating class equal to the set of cliques $\mathcal{A} = \{ab, bc, cd, da\}$, and we have

$$\mathcal{D} = \{a, b, c, d, ab, bc, cd, da\},$$

$$\mathcal{E}_\ominus = \{a, b, c, d, ab, bc, cd, da, ac, bd, abc, bcd, cda, dab, abcd\}.$$

The linear constraints on $\theta_E, E \notin \mathcal{D}$ are

$$\theta(ac) = \theta(bd) = \theta(abc) = \theta(bcd) = \theta(cda) = \theta(dab) = \theta(abcd) = 0$$

and the constraints on the cell probabilities are as follows:

$$p(ac) = \frac{p(a)p(c)}{p_\varnothing}, \qquad p(bd) = \frac{p(b)p(d)}{p_\varnothing},$$

$$p(abc) = \frac{p(ab)p(bc)}{p(b)}, \qquad p(bcd) = \frac{p(bc)p(cd)}{p(c)},$$

$$p(cda) = \frac{p(cd)p(da)}{p(d)}, \qquad p(dab) = \frac{p(da)p(ab)}{p(a)},$$

$$p(abcd) = \frac{p(ab)p(bc)p(cd)p(da)p_\varnothing}{p(a)p(b)p(c)p(d)}.$$

According to (2.20), we have

$$\begin{aligned}
p_\varnothing^{-1} = {} & 1 + e^{\theta(a)} + e^{\theta(b)} + e^{\theta(c)} + e^{\theta(d)} + e^{\theta(a)+\theta(b)+\theta(ab)} \\
& + e^{\theta(b)+\theta(c)+\theta(bc)} + e^{\theta(c)+\theta(d)+\theta(cd)} + e^{\theta(d)+\theta(a)+\theta(da)} \\
& + e^{\theta(a)+\theta(b)+\theta(c)+\theta(ab)+\theta(bc)} + e^{\theta(b)+\theta(c)+\theta(d)+\theta(bc)+\theta(cd)} \\
& + e^{\theta(c)+\theta(d)+\theta(a)+\theta(cd)+\theta(da)} + e^{\theta(d)+\theta(a)+\theta(b)+\theta(da)+\theta(ab)} \\
& + e^{\theta(a)+\theta(b)+\theta(c)+\theta(d)+\theta(ab)+\theta(bc)+\theta(cd)+\theta(da)}.
\end{aligned}$$

The other cell probabilities can be written in terms of $\theta$ according to (2.19).

**3. The conjugate prior for the log-linear parameter $\theta$.** From (2.22), it is clear that, for the three nested classes of models considered in this paper (graphical with respect to $G$ decomposable, graphical with respect to an arbitrary undirected graph $G$ and hierarchical) the conjugate prior for $\theta_\mathcal{D}$, as given by Diaconis and Ylvisaker [10], is given immediately by its density with respect to the Lebesgue measure

$$\begin{aligned}
(3.1) \quad \pi_\mathcal{D}(\theta_\mathcal{D}|s, \alpha) = {} & I_\mathcal{D}(s, \alpha)^{-1} \\
& \times \exp\Big\{ \sum_{D \in \mathcal{D}} \sum_{i_D \in \mathcal{I}_D^*} \theta(i_D) s(i_D) \\
& \qquad - \alpha \log\Big(1 + \sum_{E \in \mathcal{E}_\ominus} \sum_{i_E \in \mathcal{I}_E^*} \exp \sum_{F \subseteq \mathcal{D}E} \theta(i_F)\Big)\Big\},
\end{aligned}$$

where $I_{\mathcal{D}}(s, \alpha)$ is the corresponding normalising constant and

$$(3.2) \qquad (s, \alpha) = (s(i_D), D \in \mathcal{D}, i_D \in \mathcal{I}_D^*, \alpha), \qquad s \in \mathbb{R}^{d_{\mathcal{D}}}, \ \alpha \in \mathbb{R},$$

are the hyperparameters.

In order to be able to use this prior in practice, we need to answer a number of questions. The first basic question is to know for which values of the hyperparameters $(s, \alpha)$ the distribution is proper; that is, when does $I_{\mathcal{D}}(s, \alpha) < +\infty$ hold. Next, we can ask how to construct such hyperparameters. We address these two questions first and then we will give an example showing how to choose $(s, \alpha)$ to reflect prior knowledge.

LEMMA 3.1. *The prior distribution (3.1) with hyperparameters $(s, \alpha)$, as defined in (3.2), is proper if and only if $\frac{s}{\alpha}$ belongs to the $\mathcal{D}$-marginal cell probability space of $\mathcal{F}_{\mu_{\mathcal{D}}}$; that is, if and only if $\alpha > 0$ and there exists an array of real numbers $\rho(j) > 0, j \in \mathcal{I}$ such that*

$$(3.3) \qquad s(i_D) = \alpha \sum_{j_D = i_D} \rho(j), \qquad D \in \mathcal{D}, \ i_D \in \mathcal{I}_D^*,$$

*where, for $E \in \mathcal{E}$,*

$$(3.4) \qquad \rho(i(E)) = \frac{\exp \sum_{F \subseteq_{\mathcal{D}} E} \theta(i_F)}{1 + \sum_{E \in \mathcal{E}_{\ominus}} \sum_{j_E \in \mathcal{I}_E^*} \exp(\sum_{F \subseteq_{\mathcal{D}} E} \theta(j_F))}$$

*for some $\theta_{\mathcal{D}} \in \mathbb{R}^{d_{\mathcal{D}}}$.*

PROOF. Since the parameter space of (2.21) is $\Theta_{\mathcal{D}} = \mathbb{R}^{d_{\mathcal{D}}}$, by Theorem 1 of Diaconis and Ylvisaker [10], a necessary and sufficient condition for $I_{\mathcal{D}}(s, \alpha)$ to be finite is that $\alpha$ be a positive scalar and that $\frac{N}{\alpha} s = \frac{N}{\alpha}(s(i_D), D \in \mathcal{D}, i_D \in \mathcal{I}_D)$ be in the interior of the convex hull of the support of $\mu_{\mathcal{D}}$. Since the Laplace transform

$$L_{\mu_{\mathcal{D}}}(\theta_{\mathcal{D}}) = \left(1 + \sum_{E \in \mathcal{E}_{\ominus}} \sum_{j_E \in \mathcal{I}_E^*} \exp\left(\sum_{F \subseteq_{\mathcal{D}} E} \theta(j_F)\right)\right)^N$$

is defined on $\Theta_{\mathcal{D}}$, which is open, the interior of the convex hull of the support of $\mu_{\mathcal{D}}$ is equal to the mean space $M_{\mathcal{D}}$ of $\mathcal{F}_{\mu_{\mathcal{D}}}$. We therefore want to identify $M_{\mathcal{D}}$. Let $k_{\mu_{\mathcal{D}}}(\theta_{\mathcal{D}}) = \log L_{\mu_{\mathcal{D}}}(\theta_{\mathcal{D}})$. Since $\mathcal{F}_{\mu_{\mathcal{D}}}$ is a natural exponential family with parameter $\theta_{\mathcal{D}} \in \Theta_{\mathcal{D}}$, we have

$$(3.5) \qquad \begin{aligned} M_{\mathcal{D}} = \Big\{ m = (m(i_D), D \in \mathcal{D}, i_D \in \mathcal{I}_D^*) \Big| \\ m(i_D) = E(n(i_D)) = N \frac{dk_{\mu_{\mathcal{D}}}(\theta_{\mathcal{D}})}{d\theta(i_D)} = N \sum_{j \in \mathcal{I}, j_D = i_D} p(j) \Big\}, \end{aligned}$$

where $p(j)$ are as in (2.19). Therefore, $\frac{s(i_D)}{\alpha}, D \in \mathcal{D}, i_D \in \mathcal{I}_D^*$ must have the same properties as $\sum_{j \in \mathcal{I}, j_D = i_D} p(j)$, and the lemma follows. □

From Lemma 3.1, we know that, for $\pi_{\mathcal{D}}(\theta_{\mathcal{D}}|(s,\alpha))$ to be proper, we can choose $(s,\alpha)$ so that $\frac{s}{\alpha} = \frac{1}{\alpha}(s(i_D), D \in \mathcal{D}, i_D \in \mathcal{I}_D^*)$ as the $\mathcal{D}$-marginal probabilities of a fictive probability table with cell probabilities of the form (2.19) and (2.20) or equivalently the cell probabilities of a hierarchical log-linear model with generating set $\mathcal{D}$.

A first way to build the cell probabilities of such a fictive contingency table, that is to obtain $(s,\alpha)$, is to follow the lemma above:

1. Choose an arbitrary $\theta_{\mathcal{D}} = (\theta(i_D), D \in \mathcal{D}, i_D \in \mathcal{I}_D^*)$.
2. For each $i = i(E), E \in \mathcal{E}_\ominus, i_E \in \mathcal{I}_E^*$, define $\frac{p(i)}{\alpha}$ to be equal to the right-hand side of (2.19), and to the right-hand side of (2.20) for $E = \varnothing$. This defines $\frac{p(i)}{\alpha}$ for all $i \in \mathcal{I}$.
3. For $D \in \mathcal{D}, i_D \in \mathcal{I}_D^*$, let $\frac{s(i_D)}{\alpha} = \sum_{j \in \mathcal{I}, j_D = i_D} \frac{p(j)}{\alpha}$.
4. Choose $\alpha > 0$ arbitrarily and derive $s = (s(i_D), D \in \mathcal{D}, i_D \in \mathcal{I}_D^*)$ from step 3 above.

We note here that these hyperparameters are consistent across models, in the sense that the fictive marginal counts for different models can be obtained from a single $\theta = (\theta(i_D), D \subseteq_\ominus V, i \in \mathcal{I})$. For each model, that is, each $\mathcal{D}$, we then build the cell probabilities of a fictive table of counts through steps 1–3 above and the cell counts through step 4.

A second way to construct the cell probabilities of the fictive contingency table is to start with an arbitrary given contingency table with all cell counts $(\nu) = (\nu(i_D), D \subseteq_\ominus V, i \in \mathcal{I})$ positive, not necessarily integers. Let $\alpha$ denote the total count in that table. The maximum likelihood estimate $\hat{p}_{\mathcal{D}}$ of $p_{\mathcal{D}}$ of the fictive table cell probabilities satisfying the constraints of the model and the likelihood equations

$$\nu(i_D) = \alpha \sum_{j \in \mathcal{I}, j_D = i_D} \hat{p}(j), \qquad D \in \mathcal{D}, \ i_D \in \mathcal{I}_D^*,$$

exists; therefore, $s = (\nu(i_D), D \in \mathcal{D}, i_D \in \mathcal{I}_D^*)$ and $\alpha$ satisfy the conditions of Lemma 3.1. The hyperparameters obtained by this second method are also consistent across models in the sense that they are obtained from one single arbitrary given table of counts. As will be illustrated in the Spina Bifida example of Section 3.1, the first method can be very convenient.

The choice of $\alpha > 0$ in the methods given above is indeed arbitrary, but it is not innocent in the sense that, given a model determined by $\mathcal{D}$, the choice of $\alpha$ can change the shape of the prior distribution $\pi_{\mathcal{D}}(\theta_{\mathcal{D}}|(s,\alpha))$ and thus affect the posterior density and further inference. The following example illustrates our point. We will also see the impact of the choice of $\alpha$ in Section 5.

EXAMPLE 3.1.   Let us consider the graph $G$ which has the vertex set $V = \{a, b, c\}$ and two cliques $\{a, b\}$ and $\{b, c\}$. Let us also assume, for simplicity, that each variable $X_a, X_b, X_c$ is binary. The cell probability parameter $p_{\mathcal{D}}$, as defined in (2.23), is

$$p_{\mathcal{D}} = (p(a), p(b), p(c), p(ab), p(bc)).$$

We consider the graphical model Markov with respect to $G$. It is difficult to see the impact of the choice of $\alpha$ on $\pi_{\mathcal{D}}(\theta_{\mathcal{D}}|(s, \alpha))$, since we are not familiar with this distribution, but things are clearer when we look at the induced density on $p_{\mathcal{D}}$. As we shall see in Example 4.1, the density induced from $\pi_{\mathcal{D}}(\theta_{\mathcal{D}}|(s, \alpha))$ on $p_{\mathcal{D}}$ is equal to

$$\begin{aligned}
\pi_{\mathcal{D}}^p(p_{\mathcal{D}}|(s, \alpha)) = {}& \left( I_{\mathcal{D}}(s, \alpha) \left( 1 - \frac{p(a)p(c)}{p_{\varnothing}^2} \right) \right)^{-1} \\
& \times p(a)^{s(a)-s(ab)-1} p(b)^{s(b)-s(ab)-s(bc)-1} \\
& \times p(c)^{s(c)-s(bc)-1} p(ab)^{s(ab)-1} p(bc)^{s(bc)-1} \\
& \times p_{\varnothing}^{\alpha-s(a)-s(b)-s(c)+s(ab)+s(bc)-1}.
\end{aligned}$$

In order to obtain the hyperparameters, following the second method given above, we can take a fictive probability table with all entries equal. If we choose $\alpha = 1$, then $(s(a), s(b), s(c), s(ab), s(bc)) = \frac{1}{8}(4, 4, 4, 2, 2)$, while, if we choose $\alpha = 16$, $(s(a), s(b), s(c), s(ab), s(bc)) = 2(4, 4, 4, 2, 2)$. The corresponding conjugate priors for $p_{\mathcal{D}}$ are, respectively,

$$\begin{aligned}
\pi_{\mathcal{D}}^p & \left( p_{\mathcal{D}} \middle| \left( \left( \frac{1}{2}, \frac{1}{2}, \frac{1}{2}, \frac{1}{4}, \frac{1}{4} \right), 1 \right) \right) \\
& \propto p(a)^{-6/8} p(b)^{-1} p(c)^{-6/8} p(ab)^{-7/8} p(bc)^{-7/8} \left( 1 - \frac{p(a)p(c)}{p_{\varnothing}^2} \right)^{-1}, \\
\pi_{\mathcal{D}}^p & (p_{\mathcal{D}}|((8, 8, 8, 4, 4), 16)) \\
& \propto p(a)^3 p(b)^{-1} p(c)^3 p(ab) p(bc) \left( 1 - \frac{p(a)p(c)}{p_{\varnothing}^2} \right)^{-1}.
\end{aligned}$$

The ratio of the two densities is

$$\frac{I_{\mathcal{D}}((2(4, 4, 4, 2, 2), 16))}{I_{\mathcal{D}}(1/4(2, 2, 2, 1, 1), 1)} p(a)^{30/8} p(c)^{30/8} p(ab)^{15/8} p(bc)^{15/8}.$$

This ratio varies as $p_{\mathcal{D}}$ varies, and the two densities clearly have very different shapes and, therefore, give more prior weights to different $p_{\mathcal{D}}$.

Let us note here that for the example above, the underlying graph is decomposable and, as we shall see further in Section 4.1, in that case, the prior $\pi_{\mathcal{D}}^p(p_{\mathcal{D}}|(s,\alpha))$ coincides with the Hyper Dirichlet defined by Dawid and Lauritzen [8]. Using the notation of Section 4.1 and Proposition 4.1 and calling $C_1$ the clique $\{a,b\}$, $C_2$ the clique $\{b,c\}$ and $S$ the separator $\{b\}$, the ratio above can be written as the ratio of two hyper Dirichlet with hyperparameters equal to, for $\alpha = 1$,

$$\alpha^{C_l}(D) = \tfrac{1}{4}, \qquad D \subseteq C_l, \qquad l = 1, 2, \qquad \alpha^S(b) = \alpha^S(\varnothing) = \tfrac{1}{2},$$

and, for $\alpha = 16$, equal to

$$\alpha^{C_l}(D) = 4, \qquad D \subseteq C_l, \qquad l = 1, 2, \qquad \alpha^S(b) = \alpha^S(\varnothing) = 8.$$

The ratio is therefore equal to

$$\frac{\Gamma(16)\Gamma(8)^2\Gamma(1/4)^8}{\Gamma(1)\Gamma(1/2)^2\Gamma(4)^8}\left[\frac{p^{C_1}(ab)p^{C_1}(a)p^{C_1}(b)p^{C_1}(\varnothing)}{p^S(b)p^S(\varnothing)}\right]^{15/4}$$

$$\times \left[\frac{p^{C_2}(bc)p^{C_2}(b)p^{C_2}(c)p^{C_2}(\varnothing)}{p^S(b)p^S(\varnothing)}\right]^{15/4}.$$

Though expressed here in the more familiar marginal clique and separator marginal cell probabilities, this second expression of the ratio of prior densities may be more difficult to apprehend than the first one in terms of cell probabilities.

We note also that for the saturated model, the prior with $\alpha = 1$ and all fictive cell counts equal, is the prior advocated by Perks [25] (see also Dellaportas and Forster [9]).

3.1. *Moments of the cell probabilities.* We now compute the moments of the cell probabilities. As we show below through an example, these moments can, in some instances, be used to guide our choice of hyperparameters when we have prior information.

PROPOSITION 3.1. *Consider the distribution $\pi_{\mathcal{D}}(\theta_{\mathcal{D}}|(s,\alpha))$ as defined in (3.1). Let $r$ be a positive integer. Then, for $D \in \mathcal{D}, i_D \in \mathcal{I}_D^*$, the $r$th moment of the generalized odds ratio is*

$$E_{\pi_{\mathcal{D}}(\theta_{\mathcal{D}}|(s,\alpha))}\left(\left(\prod_{F \subseteq D} p(i(F))^{(-1)^{|D \setminus F|}}\right)^r\right) = E_{\pi_{\mathcal{D}}(\theta_{\mathcal{D}}|(s,\alpha))}(e^{r\theta(i_D)})$$

$$(3.6)$$

$$= \frac{I_{\mathcal{D}}(\tilde{s}_D, \alpha)}{I_{\mathcal{D}}(s, \alpha)},$$

*where the components of $\tilde{s}_D$ are equal to those of $s$ except for*

$$\tilde{s}_D(i(D)) = s(i(D)) + r.$$

*Moreover, for all $E \in \mathcal{E}$, the rth moment of the cell probabilities $p(i(E))$ is*

$$(3.7) \qquad E_{\pi_{\mathcal{D}}(\theta_{\mathcal{D}}|(s,\alpha))}(p(i(E))^r) = \frac{I_{\mathcal{D}}(\tilde{s}_{E,r}, \alpha + r)}{I_{\mathcal{D}}(s,\alpha)},$$

*where the components of $\tilde{s}_{E,r}$ are equal to those of $s$, except for*

$$\tilde{s}_{E,r}(i(F)) = s(i(F)) + r, \qquad if\ F \subseteq_{\mathcal{D}} E.$$

The proof of this proposition is simple and is omitted. Equation (3.7) follows from the fact that

$$(3.8) \qquad p(i(E)) = e^{\theta(i_{\varnothing}) + \sum_{F \subseteq \mathcal{D}^E} \theta(i_F)} = e^{\sum_{F \subseteq \mathcal{D}^E} \theta(i_F) - k_{\mu_{\mathcal{D}}(\theta_{\mathcal{D}})}}.$$

As we shall see in Section 4, the normalising constant $I_{\mathcal{D}}$ can be computed explicitly when the model is Markov with respect to a decomposable graph $G$. Otherwise, the normalising constants have to be computed numerically by the usual approximation methods.

We now show through an example how the results in Proposition 3.1 above can be used to guide our choice of $(s, \alpha)$ in the prior distribution. We consider the data given by Hook, Albright and Cross [16] and used by King and Brooks [17] to illustrate the fact, as we do it here, that with their prior they can translate prior information into values for the hyperparameters. In this dataset, there are three variables $a$, $b$ and $c$ each taking the values 1 or 0 representing the presence or absence of, respectively, birth certificates, death certificates and medical records for each individual. The individuals under study are children with spina bifida. The data consists of an incomplete contingency table for each one of six years. From Hook, Albright and Cross [16], it can reasonably be assumed that the model is the decomposable graphical model with cliques $\{a\}$ and $\{b, c\}$. Consultation with experts suggests that the interaction between factors $b$ and $c$ is negative and the presence of this negative interaction is expected to create a relative decrease in the $(bc)$ cell probability by a proportion in the interval $[0.1, 0.9]$. It was also expected that the total number of babies born with spina bifida during the study period would lie in the interval $[9, 56]$, and it was thought reasonable that a prior mean number of babies should lie in the interval $[29, 35]$.

Let us now express this prior information in terms of restrictions on $(s, \alpha)$. Since $\alpha$ can be thought of as the total count for a fictive prior contingency table, the belief about the total count could be immediately translated as, say, $\alpha = 30$, which lies in the interval $[29, 35]$. To reflect the negative interaction between factors $b$ and $c$, we chose $\theta_{bc}$ negative, and to reflect the fact that the negative interaction given by $\theta_{bc}$ would cause the ratio of the expected value of $p_{bc}$ under $\mathcal{D} = \{a, b, c\}$ and the expected value of $p_{bc}$ under $\mathcal{D} = \{a, b, c, bc\}$ to be in the interval $[0.1, 0.9]$, we want the values of

$\theta(a), \theta(b), \theta(c), \theta(bc)$ to satisfy

$$(3.9) \qquad 0.1 \leq \frac{E_0(p(bc))}{E_{\text{int}}(p(bc))} \leq 0.9,$$

where $E_{\text{int}}(p_E)$ denotes the expected value of $p_E$ under $\pi_{a,b,c,bc}$ and $E_0(p_E)$ denotes the expected value of $p_E$ under $\pi_{a,b,c}$. Equation (3.7) gives the two expected values in (3.9) in terms of the normalising constant of (3.1). Since the chosen model, with main effects and $bc$ interaction, is a decomposable graphical model, the normalising constant $I_{\mathcal{D}}(s, \alpha)$ can be obtained explicitly. Its formula is given further in Proposition 4.1. Using this formula, for both the model with and without the $bc$ interaction, that is with $\mathcal{D} = \{a, b, c, bc\}$ or $\mathcal{D} = \{a, b, c\}$, respectively, and also using Proposition 3.1, straightforward calculations yield

$$(3.10) \qquad \begin{aligned} E_{\text{int}}(p(bc)) &= \left(1 - \frac{s(a)}{\alpha}\right)\frac{s(bc)}{\alpha}, \\ E_0(p(bc)) &= \left(1 - \frac{s(a)}{\alpha}\right)\frac{s(b)}{\alpha}\frac{s(c)}{\alpha}. \end{aligned}$$

From Lemma 3.1 and from the first method given in Section 3.2, we know that if we generate arbitrary values of $\theta_{\mathcal{D}} = (\theta_a, \theta_b, \theta_c, \theta_{bc})$ and compute the quantities

$$(3.11) \qquad \frac{\rho(E)}{\alpha} = \frac{\exp\sum_{F \subseteq \ominus E}\theta(F)}{1 + \sum_{E \in \mathcal{E}_\ominus}\exp(\sum_{F \subseteq \mathcal{D}E}\theta(F))}$$

for each $E \in \{a, b, c, ab, ac, bc, abc\}$, then the prior (3.1) with hyperparameter $(s, \alpha)$ defined by

$$(3.12) \qquad \frac{s(D)}{\alpha} = \sum_{F \supseteq D}\rho(F), \qquad D \in \mathcal{D} = \{a, b, c, bc\}, \qquad \alpha > 0,$$

is proper. We can generate the $\theta_{\mathcal{D}}$ in any way. Here, we choose to generate $\theta_{bc}$ from a normal with mean $-1.12$ and variance $\frac{4}{9}$ and to generate $\theta_a, \theta_b$ and $\theta_c$ from independent normals with mean 0 and variance 1. We constrain those values to satisfy (3.9) expressed in terms of $\theta_a, \theta_b, \theta_c, \theta_{bc}$ using (3.10), (3.11) and (3.12). We choose, arbitrarily, the following set of values satisfying the required constraints:

$$\begin{aligned} \theta(a) &= -0.1200, & \theta(b) &= 1.1100, \\ \theta(c) &= -0.0100, & \theta(bc) &= -1.8800. \end{aligned}$$

This yields values of $s$ as follows:

$$s(a) = 9.2324, \qquad s(b) = 16.4593, \qquad s(c) = 16.4476, \qquad s(bc) = 6.9874,$$

which, together with $\alpha = 30$, defines a proper prior (3.1) reflecting our prior beliefs (3.9) and $\alpha \in [29, 35]$.

Clearly, if the prior information suggests that the model is not decomposable, the moments can no longer be obtained explicitly. For given values of $(s, \alpha)$, $I_{\mathcal{D}}(s, \alpha)$ has to be computed numerically and the formulae in Proposition 3.1 could only be used to verify that a particular choice of $(s, \alpha)$ is in accord with our prior belief. However, a first choice of $(s, \alpha)$ could be made from an approximating decomposable model.

3.2. *The strong hyper Markov property for graphical models.* Let us now assume that the multinomial distribution of the contingency cell counts is Markov with respect to an arbitrary undirected graph $G$. In this subsection, we will show that the generalized hyper Dirichlet in (3.1) is strong hyper Markov in the following sense. Let $P_1, \ldots, P_k$ a perfect sequence of the prime components of $G$, and let $S_2, \ldots, S_k$ be the corresponding separators. Though the prime components do not have to be complete, the separators $S_l = (\bigcup_{j=1}^{l-1} P_j) \cap P_l, l = 2, \ldots, k$ are complete by definition. We will use the notation

$$R_l = P_l \setminus \left( \bigcup_{j=1}^{l-1} P_j \right) = P_l \setminus S_l, \qquad l = 2, \ldots, k,$$

for the residuals, the notation

$$\mathcal{D}^{P_l}, \qquad l = 1, \ldots, k, \qquad \mathcal{D}^{S_l}, \ \mathcal{D}^{R_l}, \qquad l = 2, \ldots, k,$$

for the collection of complete subsets of the induced graphs $G_{P_l}, G_{S_l}, G_{R_l}$, respectively, and the notation

$$\begin{aligned}
(3.13) \qquad &\theta(\mathcal{D}^{P_l}), \qquad l = 1, \ldots, k, \\
&\theta(i_{S_l}, \mathcal{D}^{R_l}), \qquad i_{S_l} \in \mathcal{I}_{S_l}, \ l = 2, \ldots, k,
\end{aligned}$$

for, respectively, the log-linear parameters of the $P_l$-marginal multinomial and of the $R_l$-conditional multinomial given the value $i_{S_l}$ of the $S_l$-marginal cell. More precisely, these parameters are

$$\theta(\mathcal{D}^{P_l}) = (\theta^{P_l}(i_D), D \subseteq_{\mathcal{D}^{P_l}} P_l, i_D \in \mathcal{I}_D^*),$$

$$\theta(i_{S_l}, \mathcal{D}^{R_l}) = (\theta^{R_l | i_{S_l}}(i_D), D \subseteq_{\mathcal{D}^{R_l}} R_l, i_D \in \mathcal{I}_D^*),$$

where

$$\theta^{P_l}(i_D) = \log \prod_{F \subseteq D} (p^{P_l}(i_F, i_{P_l \setminus F}^*))^{(-1)^{|D \setminus F|}},$$

$$\theta^{R_l | i_{S_l}}(i_D) = \log \prod_{F \subseteq D} (p^{R_l | i_{S_l}}(i_F, i_{R_l \setminus F}^*))^{(-1)^{|D \setminus F|}}$$

and $p^{P_l}$ denotes $P_l$-marginal probabilities and $p^{R_l|i_{S_l}}$ denotes $R_l$-conditional probabilities given the values $i_{S_l}$ of the $S_l$-marginal cell.

We will say that (3.1) is strong hyper Markov with respect to $G$ if, under (3.1), the variables

$$\theta(\mathcal{D}^{P_1}),\ \theta(i_{S_l}, \mathcal{D}^{R_l}), \qquad i_{S_l} \in \mathcal{I}_{S_l},\ l = 2, \ldots, k$$

are mutually independent. This is clearly a generalization of the strong hyper Markov property as given by Dawid and Lauritzen [8]. We have the following result.

THEOREM 3.1. *If $\theta_{\mathcal{D}}$ follows the generalized hyper Dirichlet as defined in (3.1), then the joint distribution of the parameters in (3.13) has density*

$$\frac{\prod_{l=2}^{k} I_{\mathcal{D}^{S_l}}(s^{S_l}, \alpha) \prod_{l=1}^{k} \exp\{\langle \theta(\mathcal{D}^{P_l}), s(\mathcal{D}^{P_l}) \rangle - \alpha k(\theta(\mathcal{D}^{P_l}))\}}{\prod_{l=1}^{k} I_{\mathcal{D}^{P_l}}(s^{P_l}, \alpha) \prod_{l=2}^{k} \exp\{\langle \theta(\mathcal{D}^{S_l}), s(\mathcal{D}^{S_l}) \rangle - \alpha k(\theta(\mathcal{D}^{S_l}))\}}$$

$$= \frac{\exp\{\langle \theta(\mathcal{D}^{P_1}), s(\mathcal{D}^{P_1}) \rangle - \alpha k(\theta(\mathcal{D}^{P_1}))\}}{I_{\mathcal{D}^{P_1}}(s^{P_1}, \alpha)}$$

(3.14)
$$\times \prod_{l=2}^{k} \prod_{i_S \in \mathcal{I}_{S_l}} \frac{1}{I_{\mathcal{D}^{R_l}}(s^{R_l}, s(i_{S_l}))}$$

$$\times \exp\{\langle \theta(i_{S_l}, \mathcal{D}^{R_l}), s(i_{S_l}, \mathcal{D}^{R_l}) \rangle$$
$$- s(i_{S_l}) k(\theta(i_{S_l}, \mathcal{D}^{R_l}))\},$$

*where $s^A = (s(i_D), D \in \mathcal{D}^A, i_D \in \mathcal{I}_D^*)$ for $A = P_l, R_l$ or $S_l$, and*

$$\langle \theta(\mathcal{D}^{P_l}), s(\mathcal{D}^{P_l}) \rangle = \sum_{D \subseteq_{\mathcal{D}} P_l} \sum_{i_D \in \mathcal{I}_D^*} \theta^{P_l}(i_D) s(i_D),$$

$$\langle \theta(i_{S_l}, \mathcal{D}^{R_l}), s(i_{S_l}, \mathcal{D}^{R_l}) \rangle = \sum_{D \subseteq_{\mathcal{D}} R_l} \sum_{i_D \in \mathcal{I}_D^*} \theta^{R_l|i_{S_l}}(i_D) s(i_{S_l}, i_D),$$

$$k(\theta(\mathcal{D}^{P_l})) = \log\left(1 + \sum_{D \subseteq_{\ominus} P_l} \sum_{i_D \in \mathcal{I}_D^*} \exp \sum_{F \subseteq_{\mathcal{D}} P_l} \theta(i_F)\right),$$

$$k(\theta(i_{S_l}, \mathcal{D}^{R_l})) = \log\left(1 + \sum_{D \subseteq_{\ominus} R_l} \sum_{i_D \in \mathcal{I}_D^*} \exp \sum_{F \subseteq_{\mathcal{D}} R_l} \theta^{R_l|i_{S_l}}(i_F)\right).$$

*The parameters in (3.13) are therefore independently distributed; that is, (3.1) is strong hyper Markov.*

The proof is long and tedious but without conceptual difficulties and is ommitted here. It is important to note that $s(i_D), D \in \mathcal{D}, i_D \in \mathcal{I}_D^*$ is always the marginal count in the fictive contingency table attached to the prior, whether it occurs in the marginal distribution of $P_l$, the marginal distribution of $S_l$ or the conditional distribution of $R_l$ given $i_{S_l}$.

**4. The induced prior on the cell probabilities.** In this section, we give the expression of the induced conjugate prior in terms of the cell probability parameter, first for graphical models Markov with respect to a decomposable $G$ showing that we obtain the hyper Dirichlet, then for general hierarchical models, which includes, in particular, models Markov with respect to an arbitrary graph $G$. The proofs of all our results are given in the Appendix.

4.1. *Decomposable graphical models.* We first consider the case of the multinomial Markov, with respect to the decomposable graph $G$ with set of cliques $\mathcal{C} = \{C_l, l = 1, \ldots, k\}$ and set of minimal separators $\mathcal{S} = \{S_l, l = 2, \ldots, k\}$, so that $\mathcal{D}$ is the set of all possible subsets of $\mathcal{C}$. Since in this subsection we deal with joint cell probabilities as well as $C_l$-marginal or $S_l$-marginal probabilities, in the expression $p(i_D, i_{D^c}^*), p^{C_l}(i_D, i_{D^c}^*), p^{S_l}(i_D, i_{D^c}^*)$, it will be understood that we have $D$ as a subset of, respectively, $V, C_l$ and $S_l$, and $D^c$ as the complement of $D$ in $V, C_l$ and $S_l$. We also, temporarily, do not use the $i(D)$ notation but rather the more explicit, albeit more cumbersome, $(i_D, i_{D^c}^*)$ for $i(D)$.

Dawid and Lauritzen [8] defined the standard conjugate prior in terms of the clique and separator cell probabilities

$$
\begin{aligned}
&p^{C_l}(i_D, i_{D^c}^*), && D \in \mathcal{D}^{C_l}, \ l = 1, \ldots, k, \\
&p^{S_l}(i_D, i_{D^c}^*), && D \in \mathcal{D}^{S_l}, \ l = 2, \ldots, k, \ i_D \in \mathcal{I}_D^*,
\end{aligned}
\tag{4.1}
$$

and called it the hyper Dirichlet. Its density is equal to

$$
\frac{\prod_{l=1}^{k} \mathrm{Dir}_{C_l}(p_{\varnothing}^{C_l}, p^{C_l}(i_D, i_{D^c}^*); \alpha_{\varnothing}^{C_l}, \alpha^{C_l}(i_D, i_{D^c}^*), D \in \mathcal{D}^{C_l}, i_D \in \mathcal{I}_D^*)}{\prod_{l=2}^{k} \mathrm{Dir}_{S_l}(p_{\varnothing}^{S_l}, p^{S_l}(i_D, i_{D^c}^*); \alpha^{S_l}(i_D, i_{D^c}^*), D \in \mathcal{D}^{S_l}, i_D \in \mathcal{I}_D^*)}
\tag{4.2}
$$

with

$$
\begin{aligned}
&\mathrm{Dir}_{C_l}(p_{\varnothing}^{C_l}, p^{C_l}(i_D, i_{D^c}^*); \alpha^{C_l}(i_D, i_{D^c}^*), D \in \mathcal{D}^{C_l}, i_D \in \mathcal{I}_D^*) \\
&\qquad = \frac{\Gamma(\alpha_{\varnothing}^{C_l} + \sum_{D \in \mathcal{D}^{C_l}} \sum_{i_D \in \mathcal{I}_D^*} \alpha^{C_l}(i_D, i_{D^c}^*))}{\Gamma(\alpha_{\varnothing}^{C_l}) \prod_{D \in \mathcal{D}^{C_l}, i_D \in \mathcal{I}_D^*} \Gamma(\alpha^{C_l}(i_D, i_{D^c}^*))} (p_{\varnothing}^{C_l})^{\alpha_{\varnothing}^{C_l} - 1} \\
&\qquad\quad \times \prod_{D \in \mathcal{D}^{C_l}, i_D \in \mathcal{I}_D^*} (p^{C_l}(i_D, i_{D^c}^*))^{\alpha^{C_l}(i_D, i_{D^c}^*) - 1}
\end{aligned}
$$

with a similar expression for $\mathrm{Dir}_{S_l}$ and where the hyper parameters

$$
\begin{aligned}
&(\alpha_{\varnothing}^{C_l}, \alpha^{C_l}(i_D, i_{D^c}^*), D \in \mathcal{D}^{C_l}, i_D \in \mathcal{I}_D^*) \quad \text{and} \\
&(\alpha_{\varnothing}^{S_l}, \alpha^{S_l}(i_D, i_{D^c}^*), D \in \mathcal{D}^{S_l}, i_D \in \mathcal{I}_D^*)
\end{aligned}
\tag{4.3}
$$

are hyperconsistent in the sense that, if $S_l = C_i \cap C_j$, the marginal distributions on $S_l$ obtained from either of the clique marginal distributions on $C_i$ or $C_j$ are the same.

In this subsection, we derive the prior induced from $\pi_{\mathcal{D}}(\theta_{\mathcal{D}}|s, \alpha)$ in (3.1) by the change of variable from $\theta_{\mathcal{D}}$, as defined in (2.18), to $p^G$, as defined below in (4.4). We choose to work with $p^G$ which is the cell parametrization expressed in terms of marginal clique probabilities rather than with $p_{\mathcal{D}}$ as in (2.23) because we want to compare the induced prior on $p^G$ with the hyper Dirichlet. The probabilities in (4.1) are not all free variables. One way to choose the free marginal probabilities is as follows:

$$
\begin{aligned}
(4.4) \quad p^G = \Bigg( & p^{C_l}(i_D, i^*_{D^c}), D \in \mathcal{D}^{C_l} \setminus \bigcup_{j=2}^{k} \mathcal{D}^{S_j}, l = 1, \ldots, k, \\
& p^{S_l}(i_D, i^*_{D^c}), D \in \mathcal{D}^{S_l}, l = 2, \ldots, k, i_D \in \mathcal{I}^*_D \Bigg).
\end{aligned}
$$

The Jacobian of the change of variable $\theta_{\mathcal{D}} \mapsto p^G$ is given in the following lemma.

LEMMA 4.1. *The Jacobian of the change of variables from $\theta_{\mathcal{D}} = (\theta(i_D),$ $D \in \mathcal{D}, i_D \in \mathcal{I}^*_D)$ as given in (2.3) to $p^G$ as given in (4.4) is*

$$
(4.5) \quad \left| \frac{d\theta}{dp^G} \right|^{-1} = \frac{\prod_{l=1}^{k} p^{C_l}_{\varnothing} \prod_{D \in \mathcal{D}^{C_l}} \prod_{i_D \in \mathcal{I}^*_D} p^{C_l}(i_D, i^*_{D^c})}{\prod_{l=2}^{k} p^{S_l}_{\varnothing} \prod_{D \in \mathcal{D}^{S_l}} \prod_{i_D \in \mathcal{I}^*_D} p^{S_l}(i_D, i^*_{D^c})}.
$$

This lemma has already been stated in Leucari [20] in a slightly different form, and we give it here without proof. The following proposition says that the induced prior on $p^G$ is the hyper Dirichlet, which was also given by Leucari [20]. Here, we additionally give the correspondence between $(s, \alpha)$ and (4.3).

PROPOSITION 4.1. *When the graph $G$ is decomposable with set of cliques $(C_i, i = 1, \ldots, k)$ and sets of minimal separators $(S_i, i = 2, \ldots, k)$, the conjugate prior induced from (3.1) is identical to the hyper Dirichlet (4.2) with hyper parameters (4.3), where*

$$
\begin{aligned}
(4.6) \quad \alpha^{C_l}(i_D, i^*_{D^c}) &= \sum_{C_l \supseteq F \supseteq D} \sum_{\substack{j_F \in \mathcal{I}^*_F \\ (j_F)_D = i_D}} (-1)^{|F \setminus D|} s(j_F), \\
\alpha^{C_l}_{\varnothing} &= \alpha + \sum_{D \subseteq C_l} (-1)^{|D|} \sum_{i \in \mathcal{I}^*_D} s(i_D),
\end{aligned}
$$

$$\alpha^{S_l}(i_D, i_{D^c}^*) = \sum_{S_l \supseteq F \supseteq D} \sum_{\substack{j_F \in \mathcal{I}_F^* \\ (j_F)_D = i_D}} (-1)^{|F \setminus D|} s(j_F),$$

(4.7)

$$\alpha_{\varnothing}^{S_l} = \alpha + \sum_{D \subseteq S_l} (-1)^{|D|} \sum_{i \in \mathcal{I}_D^*} s(i_D).$$

*Moreover,*

$$(4.8) \quad I_{\mathcal{D}}(s, \alpha) = \frac{\prod_{l=1}^{k} \Gamma(\alpha_{\varnothing}^{C_l}) \prod_{D \in \mathcal{D}^{C_l}} \prod_{i_D \in \mathcal{I}_D^*} \Gamma(\alpha^{C_l}(i_D, i_{D^c}^*))}{\Gamma(\alpha) \prod_{l=2}^{k} \Gamma(\alpha_{\varnothing}^{S_l}) \prod_{D \in \mathcal{D}^{S_l}} \prod_{i_D \in \mathcal{I}_D^*} \Gamma(\alpha^{S_l}(i_D, i_{D^c}^*))}.$$

4.2. *Connection with previous work.* From Proposition 4.1 above, from the form (4.2) of the hyper Dirichlet and the form of the multinomial Markov with respect to $G$ decomposable, we see that the prior induced on $p^G$ from the DY conjugate prior (3.1) on $\theta_{\mathcal{D}}$ has the same form as the likelihood function in terms of $p^G$. We say that this induced prior on $p^G$ is standard conjugate, following the definition of Consonni and Veronese [5] who showed that, for the one-dimensional NEFs, the prior induced from the DY conjugate prior on the canonical parametrization $\theta$ onto the mean parametrization $\mu$ is standard conjugate if and only if the NEF has a quadratic variance function. Gutierrez-Pena and Smith [15] studied the case of a multivariate NEF. They first defined two parametrizations $\phi$ and $\lambda$ to be conjugate if the standard conjugate family of priors on $\lambda$ was identical to that induced from the standard conjugate family on $\phi$ by the change of variable from $\phi$ to $\lambda$. They denoted this property $\phi \smile \lambda$. They then showed, in their Theorem 1, that $\phi \smile \lambda$ if and only if the Jacobian $|\frac{d\phi}{d\lambda}|$ is proportional to the likelihood for $\lambda$. From Lemma 4.1 and their Theorem 1, we can then immediately obtain that the induced prior on $p^G$ is the hyper Dirichlet [though we cannot obtain (4.6) and (4.7)]. Gutierrez-Pena and Smith [15] showed that the characterization of Consonni and Veronese [5] could not be extended to multivariate NEF's and conjectured that, with an extended definition of conjugacy, quadratic variance functions could characterize multivariate NEF with $\theta$ and $\mu$ conjugate.

The result in Proposition 4.1 provides a counterexample to this conjecture, since it is easy to see that the variance function of the NEF Markov with respect to $G$ decomposable is not a quadratic function of $\mu_{\mathcal{D}}$. Leucari [20] had already observed this and, also, that $\theta_{\mathcal{D}} \smile \mu_{\mathcal{D}}$. We have proved here, in addition, that $\theta_{\mathcal{D}} \smile \mu_{\mathcal{D}} \smile p^G$. The parametrization $p^G$ is of course different from $p_{\mathcal{D}}$ as defined in (2.23). This distinction is important, since $p_{\mathcal{D}}$ is not a linear function of $\mu_{\mathcal{D}}$ and, as we shall see in Example 4.1, the parametrizations $\theta_{\mathcal{D}}$ is not conjugate to the parametrization $p_{\mathcal{D}}$, even in the case of a decomposable model.

4.3. *Arbitrary graphical and hierarchical models.* To obtain the conjugate prior in terms of $p_{\mathcal{D}}$, we need to compute the Jacobian $|\frac{d\theta_{\mathcal{D}}}{dp_{\mathcal{D}}}|$. Before doing so, we need to define the following quantities. For $C \in \mathcal{D}, H \in \mathcal{E}, i_C \in \mathcal{I}_C^*, j_H \in \mathcal{I}_H^*$, let

$$(4.9) \qquad F(i_C, j_H) = \begin{cases} (-1)^{|C|-1}, & \text{if } (j_H)_C = i_C, \\ 0, & \text{otherwise,} \end{cases}$$

be the entries of a $\prod_{D \in \mathcal{D}} |\mathcal{I}_D^*| \times \prod_{H \in \mathcal{E}} |\mathcal{I}_H^*|$ matrix $F$, where the rows are indexed by $i_C \in \mathcal{I}_C^*, C \in \mathcal{D}$ and the columns by $j_H \in \mathcal{I}_H^*, H \in \mathcal{E}$. We note that the definition of $F$ implies that

$$(4.10) \qquad \sum_{C \in \mathcal{D}, i_C \in \mathcal{I}_{\mathcal{D}}^*} F(i_C, j_H) = \sum_{\substack{C \in \mathcal{D}, i_C \in \mathcal{I}_{\mathcal{D}}^* \\ (j_H)_C = i_C}} F(i_C, j_H) = \sum_{C \subseteq_{\mathcal{D}} H} (-1)^{|C|-1}.$$

In the case of binary data for $\mathcal{D}$ and $\mathcal{E}$ as given in Section 2.4, the matrix $F$ is

$$(4.11)$$

$$F = \begin{pmatrix} 0 & 1 & 0 & 0 & 0 & 1 & 0 & 0 & 1 \\ 0 & 0 & 1 & 0 & 0 & 1 & 1 & 0 & 0 \\ 0 & 0 & 0 & 1 & 0 & 0 & 1 & 1 & 0 \\ 0 & 0 & 0 & 0 & 1 & 0 & 0 & 1 & 1 \\ 0 & 0 & 0 & 0 & 0 & -1 & 0 & 0 & 0 \\ 0 & 0 & 0 & 0 & 0 & 0 & -1 & 0 & 0 \\ 0 & 0 & 0 & 0 & 0 & 0 & 0 & -1 & 0 \\ 0 & 0 & 0 & 0 & 0 & 0 & 0 & 0 & -1 \end{pmatrix}$$

$$\begin{pmatrix} 1 & 0 & 1 & 0 & 1 & 1 & 1 \\ 0 & 1 & 1 & 1 & 0 & 1 & 1 \\ 1 & 0 & 1 & 1 & 1 & 0 & 1 \\ 0 & 1 & 0 & 1 & 1 & 1 & 1 \\ 0 & 0 & -1 & 0 & 0 & -1 & -1 \\ 0 & 0 & -1 & -1 & 0 & 0 & -1 \\ 0 & 0 & 0 & -1 & -1 & 0 & -1 \\ 0 & 0 & 0 & 0 & -1 & -1 & -1 \end{pmatrix}.$$

We also need the following two lemmas. Their proof is given in the Appendix.

LEMMA 4.2. *Let $G$ be a nondecomposable prime graph, and let $\mathcal{D}$ be as in (2.17). For the matrix $F$ as described in (4.9), the sum of the entries in each column $j_H, j \in \mathcal{I}_H^*, H \in \mathcal{E}_{\ominus}$ is such that*

$$(4.12) \qquad \sum_{i_C \in \mathcal{I}_C^*, C \in \mathcal{D}} F(i_C, j_H) = \sum_{C \subseteq_{\mathcal{D}} H} (-1)^{|C|-1} = 1,$$

*if and only if the subgraph induced by $H$ is decomposable and connected.*

We are now in a position to give the expression of the Jacobian for general graphical and hierarchical models. Let

$$\mathcal{U}_\ominus = \{F \in \mathcal{E}_\ominus | F \text{ is either nondecomposable or nonconnected}\}.$$

Let $\mathcal{U} = \mathcal{U}_\ominus \cup \{\varnothing\}$ and

$$(4.13) \qquad a(H) = \left( \sum_{C \subseteq_\mathcal{D} H} (-1)^{|C|-1} - 1 \right), \qquad H \in \mathcal{E}.$$

LEMMA 4.3. *The Jacobian $J(p_\mathcal{D}) = |\frac{dp_\mathcal{D}}{d\theta_\mathcal{D}}|$ of the transformation $p_\mathcal{D} \mapsto \theta_\mathcal{D}$ is equal to*

$$(4.14) \quad J(p_\mathcal{D}) = \prod_{\substack{D \in \mathcal{D} \\ i_D \in \mathcal{I}_D^\star}} p(i(D)) \left( 1 - \sum_{\substack{H \in \mathcal{E}_\ominus \\ l_H \in \mathcal{I}_H^\star}} p(l(H)) \sum_{F \subseteq_\mathcal{D} H} (-1)^{|F|-1} \right)$$

*for general hierarchical models.*

*In the particular case of graphical models, (4.14) becomes*

$$(4.15) \qquad J(p_\mathcal{D}) = \prod_{\substack{D \in \mathcal{D} \\ i_D \in \mathcal{I}_D^\star}} p(i(D)) \left( p_\varnothing - \sum_{\substack{H \in \mathcal{U}_\ominus \\ i_H \in \mathcal{I}_H^\star}} [p(i(H)) a(H)] \right).$$

EXAMPLE 4.1. For the same model as in Example 3.1, the Jacobian (4.15) for the graphical model Markov with respect to $G$ and binary data, is

$$(4.16) \qquad J = p(a)p(b)p(c)p(ab)p(bc) \left( p_\varnothing - \frac{p(a)p(c)}{p_\varnothing} \right).$$

We note, in reference to our discussion in Section 4.2 that $J$ does not have the same form as the likelihood which is proportional to

$$p(a)^{x_a} p(b)^{x_b} p(c)^{x_c} p(ab)^{x_{ab}} p(bc)^{x_{bc}} p_\varnothing^{x_\varnothing},$$

where $x_a, x_b, x_c, x_{ab}, x_{bc}$ and $x_\varnothing$ are appropriate integers. Therefore, $\theta_\mathcal{D}$ and $p_\mathcal{D}$ are not conjugate parametrizations in the sense of Gutierez-Pena and Smith [15] as defined in Section 4.2, even in the decomposable case.

We can now give the main result of this section; that is, the conjugate prior for $p_\mathcal{D}$ induced from (3.1).

THEOREM 4.1. *The conjugate prior distribution induced from (3.1) by the change of variable $\theta_\mathcal{D} \mapsto p_\mathcal{D}$ is*

$$(4.17) \qquad \pi_\mathcal{D}^p(p_\mathcal{D}|(s,\alpha)) = \frac{K(p_\mathcal{D})^{-1}}{I_G(s,\alpha)} \prod_{D \in \mathcal{D}} \prod_{i_D \in \mathcal{I}_D^\star} p(i(D))^{\alpha(i(D))-1} p_\varnothing^{\alpha_\varnothing - 1},$$

*where*

$$(4.18) \qquad \alpha(i(D)) = \sum_{F \supseteq D, F \in \mathcal{D}} \sum_{\substack{j_F \in \mathcal{I}_F^* \\ (j_F)_D = i_D}} (-1)^{F \setminus D} s(i_F),$$

$$(4.19) \qquad \alpha_\varnothing = \alpha - \sum_{D \in \mathcal{D}} \sum_{i_D \in \mathcal{I}_D^*} (-1)^{|D|} s(i_D)$$

*and where* $K(p_{\mathcal{D}}) = \dfrac{J(p_{\mathcal{D}})}{\prod_{D \in \mathcal{D}} \prod_{i_D \in \mathcal{I}_D^*} p(i(D))}$ *and* $J(p_{\mathcal{D}})$ *is as in* *(4.14)* *for general hierarchical models and* *(4.15)* *for graphical models.*

This result follows immediately from the expression of the conjugate prior (3.1) in terms of $\theta_{\mathcal{D}}$, (2.3) and Lemma 4.3.

EXAMPLE 4.2. When the graph is the four cycle with binary data as considered in Section 2.2, $\mathcal{U} = \{ac, bd, abcd, \varnothing\}$. From (4.11), (4.13) and the constraints $\theta(E) = 0$ for $E \notin \mathcal{D}$, it follows that $a(ac) = a(bd) = 1, a(abcd) = -1$. Moreover,

$$\frac{p(ac)}{p_\varnothing} = \frac{p(a)p(c)}{p_\varnothing^2}, \qquad \frac{p(bd)}{p_\varnothing} = \frac{p(b)p(d)}{p_\varnothing^2},$$

$$\frac{p(abcd)}{p_\varnothing} = \frac{p(ab)p(bc)p(cd)p(da)}{p(a)p(b)p(c)p(d)}.$$

For

$$\alpha_\varnothing = \alpha - s(a) - s(b) - s(c) - s(d) + s(ab) + s(bc) + s(cd) + s(da),$$

we have

$$\begin{aligned}
\pi(p_{\mathcal{D}}|(s, \alpha)) = {}& I_G(s, \alpha)^{-1} p(a)^{s(a) - s(da) - s(ab) - 1} p(b)^{s(b) - s(ab) - s(bc) - 1} \\
& \times p(c)^{s(c) - s(bc) - s(cd) - 1} p(d)^{s(d) - s(cd) - s(da) - 1} p(ab)^{s(ab) - 1} \\
& \times p(bc)^{s(bc) - 1} p(cd)^{s(cd) - 1} p(da)^{s(da) - 1} \\
& \times p_\varnothing^{\alpha_\varnothing - 1} \left( 1 - \frac{p_a p_c}{p_\varnothing^2} - \frac{p_b p_d}{p_\varnothing^2} + \frac{p(ab)p(bc)p(cd)p(da)}{p(a)p(b)p(c)p(d)} \right)^{-1}.
\end{aligned}$$

**5. Example. Czech Autoworkers data.** We illustrate the use of our new priors in model selection for a classical data set previously analyzed many times in the literature. We first describe our model selection procedure, and then we describe the data and the results of our model search. The C++ code for the implementation of our methods can be obtained upon request from the authors.

5.1. *Bayesian model selection.* The Bayesian paradigm to model determination involves choosing models with high posterior probability selected from a set $\mathcal{M}$ of competing models. We associate with each candidate model $m \in \mathcal{M}$ a neighborhood $\mathrm{nbd}(m) \subset \mathcal{M}$. Any two models in $m, m' \in \mathcal{M}$ are connected through a path $m = m_1, m_2, \ldots, m_k = m'$ such that $m_j \in \mathrm{nbd}(m_{j-1})$ for $j = 2, \ldots, k$. The MC$^3$ algorithm proposed by Madigan and York [23] constructs an irreducible Markov chain $m_t$, $t = 1, 2, \ldots$ with state space $\mathcal{M}$ and equilibrium distribution $\{p(m|(n)) : m \in \mathcal{M}\}$, where $(n)$ is the data in the form of a multi-way contingency table, and $p(m|(n))$ is the posterior probability of $m$. We assume that all the models are a priori equally likely; hence, $p(m|(n))$ is proportional with the marginal likelihood $p((n)|m)$.

If the chain is in state $m_t$ at time $t$, we draw a candidate model $m'$ from a uniform distribution on $\mathrm{nbd}(m_t)$. The chain moves in state $m'$ at time $t+1$; that is, $m_{t+1} = m'$ with probability

$$(5.1) \qquad \min\left\{1, \frac{p((n)|m_{t+1})/\# \, \mathrm{nbd}(m_{t+1})}{p((n)|m_t)/\# \, \mathrm{nbd}(m_t)}\right\},$$

where $\# \, \mathrm{nbd}(m)$ denotes the number of neighbors of $m$. Otherwise, the chain does not move; that is, we set $m_{t+1} = m_t$.

The marginal likelihood of a model $m$ is given by the ratio of normalizing constants

$$(5.2) \qquad p((n)|m) = I_{\mathcal{D}_m}(y + s, N + \alpha)/I_{\mathcal{D}_m}(s, \alpha),$$

where $\mathcal{D}_m$ are the possible interactions in $m$ as in (2.14).

The evaluation of the marginal likelihoods and the specification of model neighborhoods is done with respect to the particular properties of the set of candidate models considered:

1. *Hierarchical log-linear models.* We calculate the marginal likelihood in (5.2) through the Laplace approximation (see, e.g., Tierney and Kadane [28]) to the normalizing constants for the prior and posterior distribution of log-linear model parameters. The neighborhood of a hierarchical model $m$ consists of the hierarchical models obtained from $m$ by adding one of its dual generators (i.e., minimal terms not present in the model) or deleting one of its generators (i.e., maximal terms present in the model). For details, see Edwards and Havranek [14] and Dellaportas and Forster [9].

2. *Graphical log-linear models.* We evaluate the marginal likelihood in two different ways: (i) we use the Laplace approximation to the normalizing constants $I_{\mathcal{D}_m}(y + s, N + \alpha)$ and $I_{\mathcal{D}_m}(s, \alpha)$ as we did in the hierarchical case; (ii) we decompose the independence graph $G_m$ of $m$ in its sequence of prime components and separators and compute $I_{\mathcal{D}_m}$ as in (3.14). Dobra and Fienberg [12] describe efficient algorithms for generating such a decomposition.

The normalizing constants for the complete prime components and the separators (which are required to be complete) can be obtained explicitly (see Proposition 4.1). The normalizing constants for the incomplete prime components are estimated with the Laplace approximation. The neighborhood of a graphical model is defined by the graphs obtained by adding or removing one edge from $G_m$. Since each graph has the same number of neighbors, the acceptance probability (5.1) reduces to $\min\{1, \frac{p((n)|m_{t+1})}{p((n)|m_t)}\}$.

3. *Decomposable log-linear models.* In this case, the marginal likelihood can be explicitly calculated as in Proposition 4.1. The neighborhood of a decomposable model $m$ is given by those models whose independence graphs are decomposable and are obtained by adding or deleting one edge from $G_m$. Tarantola [27] provides algorithms for determining which edges can be changed in a given decomposable graph such that the resulting graph is still decomposable. The size of the neighborhoods of two decomposable graphs that differ by one edge is not necessarily the same; thus, the acceptance probability (5.1) does not simplify as it did for graphical log-linear models.

5.2. *Results.* We study the $2^6$ Czech Autoworkers table from Edwards and Havranek [14]. This cross-classfication of 1841 men gives six potential risk factors for coronary trombosis: ($a$) smoking, ($b$) strenuous mental work, ($c$) strenuous physical work, ($d$) systolic blood pressure, ($e$) ratio of beta and alpha lipoproteins and ($f$) family anamnesis of coronary heart disease.

In the absence of any prior information, we specify a proper conjugate prior for log-linear parameters through a fictive $2^6$ table with all entries equal to $\alpha/64$ for some $\alpha > 0$. All of the log-linear models are therefore constrained to have the same effective sample size $1841 + \alpha$. We remark that this approach to constructing a conjugate prior is equivalent to eliciting hyper-Dirichlet priors (see Madigan and York [22]). While the hyper-Dirichlet priors are restricted to decomposable log-linear models, the properties of our conjugate priors extend naturally to graphical and hierarchical log-linear models.

For each $\alpha \in \{0.01, 0.1, 1, 2, 3, 32, 64, 128\}$ we perform separate searches as follows: (i) a search over decomposable graphical models, (ii) a search over graphical models with marginal likelihoods estimated through decomposing the independence graph in its prime subgraphs, (iii) a search over graphical models with marginal likelihoods estimated through a single Laplace approximation and (iv) a search over hierarchical log-linear models. The results are shown in Tables 1, 2, 3 and 4. The four searches are labeled, respectively, "Dec.," "Graph./PM," "Graph./Lapl" and "Hierar." For each search type and each value of $\alpha$, we run four separate Markov chains from a random starting model for 25,000 iterations with a burn-in of 5000 iterations. We give the models whose normalized posterior probabilities are greater than 0.05 as well as the median log-linear models that are labeled with "med." A median model contains those interaction terms having a posterior inclusion probability greater than 0.5.

TABLE 1
*The most probable log-linear models for $\alpha \in \{0.01, 0.1\}$*

| Search | $\alpha = 0.01$ | | $\alpha = 0.1$ | |
|---|---|---|---|---|
| Dec. | $ac\|bc\|d\|be\|f$ | 0.278 | $ac\|bc\|d\|be\|f$ | 0.172 |
| | $ac\|bc\|d\|e\|f$ | 0.236 | $ac\|bc\|be\|de\|f$ | 0.156 |
| | $ac\|bc\|d\|ae\|f$ | 0.212 | $ac\|bc\|d\|ae\|f$ | 0.131 |
| | $ac\|bc\|d\|ce\|f$ | 0.147 | $ac\|bc\|ae\|de\|f$ | 0.119 |
| | $ac\|bc\|d\|e\|f$ | med | $ac\|bc\|d\|e\|f$ | med |
| Graph./PM | $ac\|bc\|d\|ae\|be\|f$ | 0.856 | $ac\|bc\|d\|ae\|be\|f$ | 0.380 |
| | $ac\|bc\|ae\|be\|de\|f$ | 0.078 | $ac\|bc\|ae\|be\|de\|f$ | 0.344 |
| | | | $ac\|bc\|ad\|ae\|be\|f$ | 0.140 |
| | $ac\|bc\|d\|ae\|be\|f$ | med | $ac\|bc\|d\|ae\|be\|f$ | med |
| Graph./Lapl | $ac\|bc\|ae\|be\|de\|f$ | 0.393 | $ac\|bc\|be\|de\|f$ | 0.450 |
| | $ac\|bc\|d\|ae\|be\|f$ | 0.336 | $ac\|bc\|ad\|ae\|be\|f$ | 0.184 |
| | $ac\|bc\|ad\|ae\|be\|f$ | 0.160 | $ac\|bc\|d\|ae\|be\|f$ | 0.122 |
| | | | $ac\|bc\|be\|ade\|f$ | 0.067 |
| | $ac\|bc\|d\|ae\|be\|f$ | med | $ac\|bc\|ae\|be\|de\|f$ | med |
| Hierar. | $ac\|bc\|ad\|ae\|ce\|de\|f$ | 0.251 | $ac\|bc\|ad\|ae\|ce\|de\|f$ | 0.362 |
| | $ac\|bc\|ad\|ae\|be\|de\|f$ | 0.157 | $ac\|bc\|ad\|ae\|be\|de\|f$ | 0.227 |
| | $ac\|bc\|ae\|ce\|de\|f$ | 0.136 | $ac\|bc\|ae\|ce\|de\|f$ | 0.062 |
| | $ac\|bc\|d\|ae\|ce\|f$ | 0.116 | | |
| | $ac\|bc\|ae\|be\|de\|f$ | 0.085 | | |
| | $ac\|bc\|d\|ae\|be\|f$ | 0.0725 | | |
| | $ac\|bc\|ad\|ae\|ce\|f$ | 0.055 | | |
| | $ac\|bc\|ad\|ae\|ce\|de\|f$ | med | $ac\|bc\|ad\|ae\|ce\|de\|f$ | med |

We compare our highest posterior probability models with the log-linear models identified by Dellaportas and Forster [9], who proposed a reversible jump Markov chain Monte Carlo with normal priors for log-linear parameters, and with the decomposable models selected by Madigan and Raftery [21], who employed a hyper-Dirichlet prior for cell probabilities. Our most probable decomposable model $bc|ace|ade|f$ for $\alpha = 1, 2$ or 3 is the same decomposable model as the one identified by both Dellaportas and Forster [9] and Madigan and Raftery [21]. Our most probable graphical model $ac|bc|be|ade|f$ for $\alpha = 1, 2$ or 3 in the "Graph./Lapl" search is precisely the most probable model of Dellaportas and Forster [9] and is the second best model selected by Edwards and Havranek [14]. Similarly, our most probable hierarchical model $ac|bc|ad|ae|ce|de|f$ for $\alpha = 1, 2$ or 3 coincides with the model with the largest posterior probability identified by [9]. The same consistency of the results obtained holds for most of the highest probable models selected by us and by Dellaportas and Forster [9].

TABLE 2
*The most probable log-linear models for $\alpha \in \{1, 2\}$*

| Search | $\alpha = 1$ | | $\alpha = 2$ | |
|---|---|---|---|---|
| Dec. | $bc|ace|ade|f$ | 0.250 | $bc|ace|ade|f$ | 0.261 |
| | $bc|ace|de|f$ | 0.104 | $bc|ace|de|f$ | 0.177 |
| | $bc|ad|ace|f$ | 0.102 | $bc|ace|de|bf$ | 0.096 |
| | $ac|bc|be|de|f$ | 0.060 | $bc|ad|ace|f$ | 0.072 |
| | $bc|ace|de|bf$ | 0.051 | $bc|ace|de|bf$ | 0.065 |
| | $bc|ace|de|f$ | med | $bc|ad|ace|de|f$ | med |
| Graph./PM | $ac|bc|ae|be|de|f$ | 0.446 | $ac|bc|ae|be|de|f$ | 0.371 |
| | $ac|bc|ad|ae|be|f$ | 0.182 | $ac|bc|ad|ae|be|f$ | 0.151 |
| | $ac|bc|ae|be|de|bf$ | 0.092 | $ac|bc|ae|be|de|bf$ | 0.136 |
| | $ac|bc|d|ae|be|f$ | 0.054 | $ac|bc|ae|be|de|e|f$ | 0.057 |
| | | | $ac|bc|ad|ae|be|bf$ | 0.055 |
| | $ac|bc|ae|be|de|f$ | med | $ac|bc|ae|be|de|f$ | med |
| Graph./Lapl | $ac|bc|be|ade|f$ | 0.301 | $ac|bc|be|ade|f$ | 0.341 |
| | $ac|bc|ae|be|de|f$ | 0.203 | $ac|bc|be|ade|bf$ | 0.141 |
| | $ac|bc|be|ade|bf$ | 0.087 | $ac|bc|ae|be|de|f$ | 0.116 |
| | $ac|bc|ad|ae|be|f$ | 0.083 | $ac|bc|be|ade|e|f$ | 0.059 |
| | $ac|bc|ae|be|de|bf$ | 0.059 | | |
| | $ac|bc|ad|ae|be|de|f$ | med | $ac|bc|be|ade|f$ | med |
| Hierar. | $ac|bc|ad|ae|ce|de|f$ | 0.241 | $ac|bc|ad|ae|ce|de|f$ | 0.175 |
| | $ac|bc|ad|ae|be|de|f$ | 0.151 | $ac|bc|ad|ae|be|de|f$ | 0.110 |
| | $ac|bc|ad|ae|be|ce|de|f$ | 0.076 | $ac|bc|ad|ae|be|ce|de|f$ | 0.078 |
| | $ac|bc|ad|ae|ce|de|bf$ | 0.070 | $ac|bc|ad|ae|ce|de|bf$ | 0.072 |
| | $ac|bc|ad|ae|ce|de|f$ | med | $ac|bc|ad|ae|be|ce|de|f$ | med |

Tables 1, 2, 3 and 4 show the sensitivity of the highest posterior probability models with respect to the choice of priors and to the class of log-linear models considered. For a fixed $\alpha$, the highest probable models become sparser as we sequentially relax the structural constraints from decomposable to graphical and hierarhical. We remark that the most probable graphical (hierarchical) models can be obtained from the most probable decomposable (graphical) models by dropping some of the second-order interaction terms. Increasing $\alpha$ from 1 to 128 (i.e., increasing each fictive cell count from 1/64 to 2) leads to the inclusion of additional terms in the highest probable models. We also remark that the two estimation methods for the marginal likelihoods of graphical models yield consistent results. For $\alpha = 0.01$ and $\alpha = 0.1$, we note that our results, though not abherent, are not as entirely consistent with the results obtained for $\alpha \in \{1, 2, 3, 64, 128\}$ as these results are between themselves or with results obtained in previous studies. This

Table 3
*The most probable log-linear models for $\alpha \in \{3, 32\}$*

| Search | $\alpha = 3$ | | $\alpha = 32$ | |
|---|---|---|---|---|
| Dec. | $bc\|ace\|ade\|f$ | 0.312 | $bc\|ace\|ade\|bf$ | 0.136 |
| | $bc\|ace\|ade\|bf$ | 0.155 | $ace\|bce\|ade\|bf$ | 0.098 |
| | $bc\|ace\|de\|f$ | 0.107 | $bc\|ace\|ade\|f$ | 0.062 |
| | $bc\|ace\|ade\|ef$ | 0.065 | $abc\|ace\|ade\|f$ | 0.060 |
| | $bc\|ace\|de\|bf$ | 0.053 | $bc\|ace\|ade\|ef$ | 0.057 |
| | $bc\|ace\|ade\|f$ | med | $bc\|ace\|ade\|f$ | med |
| Graph./PM | $ac\|bc\|ae\|be\|de\|f$ | 0.316 | $bc\|ace\|ade\|f$ | 0.068 |
| | $ac\|bc\|ae\|be\|de\|bf$ | 0.157 | $ac\|bc\|ade\|bde\|bf$ | 0.050 |
| | $ac\|bc\|ad\|ae\|be\|f$ | 0.128 | | |
| | $ac\|bc\|ae\|be\|de\|ef$ | 0.066 | | |
| | $ac\|bc\|ad\|ae\|de\|bf$ | 0.064 | | |
| | $ac\|bc\|ae\|be\|de\|f$ | med | $bc\|be\|ace\|ade\|f$ | med |
| Graph./Lapl | $ac\|bc\|be\|ade\|f$ | 0.339 | $ac\|bc\|be\|ade\|bf$ | 0.188 |
| | $ac\|bc\|be\|ade\|bf$ | 0.172 | $ac\|bc\|be\|ade\|f$ | 0.103 |
| | $ac\|bc\|ae\|be\|de\|f$ | 0.077 | $ac\|bc\|be\|ade\|ef$ | 0.079 |
| | $ac\|bc\|be\|ade\|ef$ | 0.072 | $ac\|bc\|be\|ade\|af\|bf$ | 0.057 |
| | | | $ac\|bc\|be\|ade\|bf\|df$ | 0.052 |
| | $ac\|bc\|be\|ade\|f$ | med | $ac\|bc\|be\|ade\|bf$ | med |
| Hierar. | $ac\|bc\|ad\|ae\|ce\|de\|f$ | 0.137 | $ac\|bc\|ad\|ae\|be\|ce\|de\|bf$ | 0.020 |
| | $ac\|bc\|ad\|ae\|be\|de\|f$ | 0.086 | | |
| | $ac\|bc\|ad\|ae\|be\|ce\|de\|f$ | 0.075 | | |
| | $ac\|bc\|ad\|ae\|ce\|de\|bf$ | 0.069 | | |
| | $ac\|bc\|ad\|ae\|be\|ce\|de\|f$ | med | $ac\|bc\|ad\|ae\|be\|ce\|de\|bf$ | med |

is not surprising, since for values of $\alpha$ very close to 0, we encounter two potential problems: first, the unknown behaviour of the Bayes factor as $\alpha$ tends to 0 and, second, the evaluation of the prior normalizing constant. The first problem is a very important general problem regarding the behaviour of Bayes factors when the total "fictive cell counts" tends to 0 (see Steck and Jaakkola [26] for related results on directed acyclic graphs). This problem needs careful study and will be the subject of further work for our particular prior distributions. In this paper, we confine ourselves to observing a difference in the behaviour of the model selection results between the case $\alpha \in \{0.01, 0.1\}$ and the results obtained for other larger values of $\alpha$. Regarding the second problem, we have observed numerically that, for $\alpha$ close to 0, the prior distribution for each $\theta_D$ is very flat; hence, the Laplace approximation is bound to yield poor results.

A very interesting question is whether there exists evidence of a link between $f$ and the other five risk factors. Whittaker [30], page 263, chooses the graphical model $abce|ade|bf$ that links $f$ with $b$, strenuous mental work. The most probable models identified by Edwards and Havranek [14], Madigan and Raftery [21] or Dellaportas and Forster [9] indicate the independence of $f$ from the other risk factors. Their findings are consistent with the models we identify for smaller values of $\alpha$. However, as we increase the grand total $\alpha$ in the prior fictive table we use, a direct link between $b$ and $f$ appears in our highest probable models. Table 5 shows the posterior inclusion probability of the first-order interaction between $b$ and $f$ for various choices of $\alpha$ and structural model constraints. Table 5 seems to confirm Whittaker's findings, because the posterior probability of the edge $bf$ increases from 0.002 to almost 1. This edge does not appear in sparser models corresponding to smaller values of $\alpha$, because there are stronger associations among $a$, $b$, $c$, $d$ and $e$ than between $f$ and another risk factor.

TABLE 4
*The most probable log-linear models for $\alpha \in \{64, 128\}$*

| Search | $\alpha = 64$ | | $\alpha = 128$ | |
|---|---|---|---|---|
| Dec. | $ace\|bce\|ade\|bcf$ | 0.134 | $ace\|bce\|ade\|bcf$ | 0.359 |
| | $ace\|bce\|ade\|bf$ | 0.118 | $ace\|ade\|bcf\|cef$ | 0.133 |
| | $ace\|ade\|bcf$ | 0.081 | $abc\|ace\|ade\|bcf$ | 0.105 |
| | $bc\|ace\|ade\|bf$ | 0.071 | $abc\|abe\|ade\|bcf$ | 0.104 |
| | $abc\|ace\|ade\|acf$ | 0.062 | $ace\|ade\|acf$ | 0.089 |
| | $abc\|ace\|ade\|bf$ | 0.055 | $abce\|ade\|acf$ | 0.060 |
| | $abc\|abe\|ade\|acf$ | 0.052 | $ace\|ade\|bcef$ | 0.051 |
| | | | $ace\|ade\|abcf$ | 0.050 |
| | $bc\|be\|ace\|ade\|bf$ | med | $be\|ace\|ade\|bcf$ | med |
| Graph./PM | $ace\|bce\|ade\|bcf$ | 0.091 | $ace\|bce\|ade\|bcf$ | 0.280 |
| | $ace\|bce\|ade\|bf$ | 0.080 | $ace\|bce\|ade\|bde\|bcf$ | 0.138 |
| | $ace\|ade\|acf$ | 0.055 | $ace\|ade\|bcf\|cef$ | 0.104 |
| | | | $abc\|ace\|ade\|bcf$ | 0.082 |
| | | | $abc\|abe\|ade\|bcf$ | 0.081 |
| | | | $ace\|ade\|bcf$ | 0.070 |
| | $bc\|be\|ace\|ade\|bf$ | med | $be\|ace\|ade\|bcf$ | med |
| Graph./Lapl | $ac\|bc\|be\|ade\|bf$ | 0.162 | $ac\|be\|ade\|bcf$ | 0.161 |
| | $ac\|be\|ade\|bcf$ | 0.128 | $ace\|bce\|ade\|bcf$ | 0.114 |
| | $ac\|bc\|be\|ade\|af\|bf$ | 0.068 | $ac\|be\|ade\|bcf\|df$ | 0.109 |
| | $ac\|bc\|be\|ade\|bf\|df$ | 0.068 | $ace\|bce\|ade\|bcf\|df$ | 0.077 |
| | $ac\|bc\|be\|ade\|ef$ | 0.068 | $ac\|ade\|bcf\|bef$ | 0.069 |
| | $ac\|bc\|be\|ade\|f$ | 0.057 | $ac\|ade\|bcf\|bef\|def$ | 0.064 |
| | $ac\|be\|ade\|bcf\|df$ | 0.054 | | |
| | $ac\|bc\|be\|ade\|bf$ | med | $ac\|be\|ade\|bcf$ | med |
| Hierar. | $ac\|bc\|ad\|ae\|be\|ce\|de\|bf$ | 0.023 | $ac\|ad\|ae\|be\|ce\|de\|bcf\|ef$ | 0.012 |
| | $ac\|bc\|ad\|ae\|be\|ce\|de\|bf$ | med | $ac\|ad\|ae\|be\|ce\|de\|bcf\|df\|ef$ | med |

TABLE 5
*Posterior inclusion probabilities for edge $bf$*

| | $\alpha$ | | | | | | | |
|---|---|---|---|---|---|---|---|---|
| Search | **0.01** | **0.1** | **1** | **2** | **3** | **32** | **64** | **128** |
| Dec. | 0.002 | 0.022 | 0.149 | 0.219 | 0.261 | 0.49 | 0.645 | 0.918 |
| Graph./PM | 0.002 | 0.033 | 0.152 | 0.222 | 0.263 | 0.476 | 0.608 | 0.898 |
| Graph./Lapl | 0.027 | 0.080 | 0.205 | 0.260 | 0.296 | 0.570 | 0.716 | 0.940 |
| Hierar. | 0.028 | 0.084 | 0.227 | 0.297 | 0.343 | 0.708 | 0.867 | 0.995 |

**6. Further developments.** The family of conjugate priors introduced in this paper has a large area of applicability. Dobra and Massam [13] make use of these priors to analyze eight and 16-dimensional contingency tables. Due to the inherent sparsity of such datasets, penalizing for increased model complexity is key and can be done naturally in the framework we have described. The same priors are used in Dobra et al. [11] to develop variable selection approaches for regressions induced by log-linear models. Their examples involve data from genomewide studies.

## APPENDIX

**A.1. Proof of Lemma 4.1.** Let $\mu_{\mathcal{D}} = \{\mu(i_D) = E(y(i_D)), D \in \mathcal{D}, i_D \in \mathcal{I}_D^*\}$ denote the mean parameter of (2.21). Leucari (2004) has proved that the Jacobian $\left|\frac{d\theta_{\mathcal{D}}}{d\mu_{\mathcal{D}}}\right|$ can be expressed as the inverse of the right-hand side of (4.5) for the binary case. The proof can immediately be extended to the discrete data case. To complete the proof of Lemma 4.1, it remains to show that the Jacobian of $\mu_{\mathcal{D}} \mapsto p_G$ is equal to 1. This is immediate, since the parameters in (4.4) are such that

$$p^{C_l}(i_D, i_{D^c}^*) = \sum_{F \subseteq C_l \setminus D} \sum_{j_F \in \mathcal{I}_F^*} (-1)^{|F|} \mu(i_D, j_F),$$

$$p^{S_l}(i_D, i_{D^c}^*) = \sum_{F \subseteq S_l \setminus D} \sum_{j_F \in \mathcal{I}_F^*} (-1)^{|F|} \mu(i_D, j_F).$$

The Jacobian of $\mu_{\mathcal{D}} \mapsto p_G$ is therefore, 1 and the lemma is proved.

**A.2. Proof of Proposition 4.1.** The distribution of $Y$ in (2.21) is Markov with respect to $G$; therefore,

(A.1) $$p(i) = \frac{\prod_{l=1}^{k} p^{C_l}(i_{C_l})}{\prod_{l=2}^{k} p^{S_l}(i_{S_l})}.$$

Since we are dealing with joint as well as marginal probabilities, we revert to the more precise notation $(i_F, i_{F^c}^*)$ rather than $i(F)$. Then,

$$\theta(i_E) = \sum_{F \subseteq E} (-1)^{|E \setminus F|} \log p(i_F, i_{F^c}^*)$$

$$= \sum_{F \subseteq E} (-1)^{|E \setminus F|} \left( \sum_{l=1}^{k} \log p^{C_l}(i_{F \cap C_l}, i_{F^c \cap C_l}^*) - \sum_{l=2}^{k} \log p^{S_l}(i_{F \cap S_l}, i_{F^c \cap S_l}^*) \right)$$

$$= \sum_{l=1}^{k} \left( \sum_{F \subseteq E} (-1)^{|E \setminus F|} \log p^{C_l}(i_{F \cap C_l}, i_{F^c \cap C_l}^*) \right)$$

$$\quad - \sum_{l=2}^{k} \left( \sum_{F \subseteq E} (-1)^{|E \setminus F|} \log p^{S_l}(i_{F \cap S_l}, i_{F^c \cap S_l}^*) \right)$$

$$= \sum_{l=1}^{k} \theta^{C_l}(i_{E \cap C_l}) - \sum_{l=2}^{k} \theta^{S_l}(i_{E \cap S_l}).$$

If $E \subseteq C_l$, $E \cap C_l = E$ and $\theta^{C_l}(i_{E \cap C_l}) = \theta^{C_l}(i_E)$. If $E \not\subseteq C_l$, then, by Lemma 2.1, $\theta^{C_l}(i_{E \cap C_l}) = 0$ and, similarly, for $\theta^{S_l}(i_{E \cap S_l})$. We therefore have

$$(A.2) \qquad \theta(i_D) = \sum_{l=1}^{k} \theta^{C_l}(i_D) - \sum_{l=2}^{k} \theta^{S_l}(i_D), \qquad D \in \mathcal{D},$$

where

$$\theta^{C_i}(i_D) = \sum_{F \subseteq D} (-1)^{|D \setminus F|} \log p^{C_i}(i_F, i_{F^c \cap C_i}^*) \qquad \text{for } D \subseteq C_i,$$

$$\theta_{\varnothing}^{C_i} = \log p_{\varnothing}^{C_i},$$

$$\theta^{C_i}(i_D) = 0 \qquad \text{for } D \not\subseteq C_i,$$

and similar expressions for $\theta^{S_i}(i_D)$ (see also Consonni and Leucari [4] for the derivation of these formulas in the case of binary data). For the remainder of this proof, it will be understood that, for $E \subseteq C_l$, we use the notation $p^{C_l}(i_E, i_{E^c}^*)$ for $p^{C_l}(i_E, i_{E^c \cap C_l}^*)$. Now, from (A.1), we also have

$$\log p_{\varnothing} = \sum_{l=1}^{k} \log p_{\varnothing}^{C_l} - \sum_{l=2}^{k} \log p_{\varnothing}^{S_l}.$$

Therefore, (3.1) can be written as

$$\pi_G(\theta_{\mathcal{D}}(p^G)|s, \alpha)$$

$$\propto \left( \prod_{l=1}^{k} \exp \left\{ \sum_{D \in \mathcal{D}^{C_l}} \sum_{i_D \in \mathcal{I}_D^*} \left( \sum_{E \subseteq D} (-1)^{|D \setminus E|} \log p^{C_l}(i_E, i_{E^c}^*) \right) s(i_D) \right. \right.$$

$$+ \alpha \log p_\varnothing^{C_l} \Big\} \Big)$$

$$\times \left( \prod_{l=2}^{k} \exp \Big\{ \sum_{D \in \mathcal{D}^{S_l}} \sum_{i_D \in \mathcal{I}_D^*} \Big( \sum_{E \subseteq D} (-1)^{|D \setminus E|} \log p^{S_l}(i_E, i_{E^c}^*) \Big) s(i_D) \right.$$

(A.3)
$$\left. + \alpha \log p_\varnothing^{S_l} \Big\} \right)^{-1}$$

$$= \frac{\prod_{l=1}^{k} \exp \{ \sum_{E \in \mathcal{D}^{C_l}} \sum_{i_E \in \mathcal{I}_E^*} \alpha^{C_l}(i_E, i_{E^c}^*) \log p^{C_l}(i_E, i_{E^c}^*) + \alpha_\varnothing^{C_l} \log p_\varnothing^{C_l} \}}{\prod_{l=2}^{k} \exp \{ \sum_{E \in \mathcal{D}^{S_l}} \sum_{i_E \in \mathcal{I}_E^*} \alpha^{S_l}(i_E, i_{E^c}^*) \log p^{S_l}(i_E, i_{E^c}^*) + \alpha_\varnothing^{S_l} \log p_\varnothing^{S_l} \}}$$

$$= \frac{\prod_{l=1}^{k} (p_\varnothing^{C_l})^{\alpha_\varnothing^{C_l}} \prod_{E \in \mathcal{D}^{C_l}} \prod_{i_E \in \mathcal{I}_E^*} (p^{C_l}(i_E, i_{E^c}^*))^{\alpha^{C_l}(i_E, i_{E^c}^*)}}{\prod_{l=2}^{k} (p_\varnothing^{S_l})^{\alpha_\varnothing^{S_l}} \prod_{E \in \mathcal{D}^{S_l}} \prod_{i_E \in \mathcal{I}_E^*} (p^{S_l}(i_E, i_{E^c}^*))^{\alpha^{S_l}(i_E, i_{E^c}^*)}},$$

where $\alpha^{C_l}(i_E, i_{E^c}^*), \alpha^{S_l}(i_E, i_{E^c}^*), \alpha_\varnothing^{C_l}$ and $\alpha_\varnothing^{S_l}$ are as defined in (4.6) and (4.7).

The induced prior on $p^G$ is obtained by multiplying (A.3) by the Jacobian (4.5) and it follows immediately that it is the hyper Dirichlet with hyper parameters as given in (4.6) and (4.7).

The expression of (4.8) is obtained by noticing that, for any $C_l$,

$$\alpha_\varnothing^{C_l} + \sum_{E \in \mathcal{D}^{C_l}} \sum_{i_E \in \mathcal{I}_E^*} \alpha^{C_l}(i_E, i_{E^c}^*) = \alpha = \alpha_\varnothing^{S_l} + \sum_{E \in \mathcal{D}^{S_l}} \sum_{i_E \in \mathcal{I}_E^*} \alpha^{S_l}(i_E, i_{E^c}^*).$$

This completes the proof.

**A.3. Proof of Lemma 4.2.** For ease of notation, we will give the proof of the lemma in the case of binary data. Given definition (4.9), the proof for binary and discrete data are exact parallel of each other. Since, for binary data, for each $C \in \mathcal{D}$ and $H \in \mathcal{E}_\ominus$ there is only one cell in $\mathcal{I}_C^*$ and $\mathcal{I}_H^*$, respectively, we will adopt the notation

$$F_{C,H} = F(i_C, j_H), \qquad C \in \mathcal{D}, \ H \in \mathcal{E}.$$

Let us first prove that, if $H$ is decomposable, (4.12) is true. We proceed by induction on the number $k$ of cliques of $H$. Let $\mathcal{C} = \{C_1, \ldots, C_k\}$ be a perfect ordering of the cliques of $H$.

If $H$ is complete, that is, $k = 1$, we consider two cases, the case where $|H|$ is even and the case where it is odd. For $|H| = 2p, p \in \mathbf{N}$, there are $n_e = \sum_{k=1}^{p} \binom{|H|}{2k}$ nonempty subsets of $H$ of even cardinality and $n_o = \sum_{k=0}^{p-1} \binom{|H|}{2k+1}$ subsets of odd cardinality. Therefore,

$$\sum_{C \subseteq_{\mathcal{D}} H} (-1)^{|C|-1} = \sum_{k=1}^{2p} \binom{2p}{k} (-1)^{k+1} = (1-1)^{2p} - \binom{2p}{0} (-1)^1 = 1$$

and (4.12) is verified. We omit the proof for the case $|H| = 2p - 1$, which is parallel to that of the previous case. Therefore, (4.12) is verified for $k = 1$.

Let us now assume that $H$ is decomposable but not complete, that is, $k > 1$, and let us assume that (4.12) is true for any decomposable subset with $k-1$ cliques. It is well known from the theory of decomposable graphs that, if we write $H_{k-1} = \bigcup_{j=1}^{k-1} C_j$, then $H = H_{k-1} \cup (C_k \setminus S_k)$, where $S_k = H_{k-1} \cap C_k$ is the $k$th minimal separator in $H$. Therefore, we have

$$\text{(A.4)} \qquad \sum_{C \subseteq_\mathcal{D} H} F_{C,H} = \sum_{C \subseteq_\mathcal{D} H_{k-1}} F_{C,H} + \left( \sum_{C \subseteq_\mathcal{D} C_k} F_{C,H} - \sum_{C \subseteq_\mathcal{D} S_k} F_{C,H} \right).$$

The first term on the right-hand side of (A.4) is equal to 1 by our induction assumption, while each one of the two other terms is also equal to 1, because both $C_k$ and $S_k$ are complete; therefore, (4.12) is also verified for decomposable $H$.

Let us now prove that if $H$ is not decomposable and connected, $\sum_{C \subseteq_\mathcal{D} H} F_{C,H}$ cannot be equal to 1. If $H$ is not connected and its connected components $H^{(1)}, \ldots, H^{(l)}$, for some $l \geq 2$, are all decomposable, we clearly have

$$\sum_{C \subseteq_\mathcal{D} H} F_{C,H} = \sum_{j=1}^{l} \left( \sum_{C \subseteq_\mathcal{D} H^{(j)}} F_{C,H^{(j)}} \right) \neq 1.$$

If $H$ is not connected and its components are not all decomposable, this implies that there is a nondecomposable subset $F_1$ of $G$, which can be separated from another subset $F_2$ of $G$, but this contradicts our assumption that $G$ is a prime component of $G$. So, this case does not occur.

If $H$ is not decomposable and connected, consider its set of cliques $\{C_1, \ldots, C_k\}$. Since $H$ is not decomposable, there is no perfect ordering of the cliques; therefore, for any given ordering, there exists a nonempty subset $\mathcal{Q} \subseteq \{3, \ldots, k\}$ such that, for $j \in \mathcal{Q}$, there is no $i < j$ in the given ordering of the cliques of $H$ with $S_j = C_j \cap (\bigcup_{l=1}^{j-1} C_l) \subseteq C_i$. Therefore,

$$S_j = C_j \cap \left( \bigcup_{l=1}^{j-1} C_l \right) = \bigoplus_{l=1}^{s_j} S_{j_l}, \qquad 2 \leq s_j \leq j - 1,$$

where the $S_{j_l}$ can be chosen to be disjoints, with $S_{j_l} \subseteq C_j \cap C_m$ for some $m \in \{1, \ldots, j-1\}$.

For $j \in \overline{\mathcal{Q}} = \{2, \ldots, k\} \setminus \mathcal{Q}$, there exists $i < j$ in the given ordering of the cliques of $H$ such that $S_j \subseteq C_i$. Therefore,

$$\text{(A.5)} \qquad \sum_{C \subseteq_\mathcal{D} H} F_{C,H} = \sum_{C \subseteq_\mathcal{D} C_1} F_{C,H} + \sum_{j \in \overline{\mathcal{Q}}} \left( \sum_{C \subseteq_\mathcal{D} C_j} F_{C,H} - \sum_{C \subseteq_\mathcal{D} S_j} F_{C,H} \right)$$

$$\text{(A.6)} \qquad \qquad + \sum_{j \in \mathcal{Q}} \left( \sum_{C \subseteq_\mathcal{D} C_j} F_{C,H} - \sum_{l=1}^{s_j} \sum_{C \subseteq_\mathcal{D} S_{j_l}} F_{C,H} \right).$$

The sums $\sum_{C \subseteq_{\mathcal{D}} U} F_{C,H}, U = C_1, C_j, S_j, j \in \overline{\mathcal{Q}}$ are all equal to 1, since each of $C_1, C_j, S_j, j \in \overline{\mathcal{Q}}$ are complete and connected; therefore, the right-hand side of (A.5) is equal to 1. For the same reason, on line (A.6), for $U = C_j, S_{j_l}, j \in \mathcal{Q}, l = 1, \ldots, s_j, \sum_{C \subseteq_{\mathcal{D}} U} F_{C,H} = 1$. Since $s_j \geq 2$,

$$\sum_{C \subseteq_{\mathcal{D}} C_j} F_{C,H} - \sum_{l=1}^{s_j} \sum_{C \subseteq_{\mathcal{D}} S_{j_l}} F_{C,H} \leq -1, \qquad j \in \mathcal{Q}.$$

Therefore, the sum on line (A.6) is less than or equal to $-|\mathcal{Q}|$. It follows that

$$\sum_{C \subseteq_{\mathcal{D}} H} F_{C,H} \leq 0$$

and in particular it cannot be equal to 1. The lemma is now proved.

**A.4. Proof of Lemma 4.3.** The rows and columns of the Jacobian of the change of variables $p_{\mathcal{D}} \mapsto \theta_{\mathcal{D}}$ is a $d_{\mathcal{D}} \times d_{\mathcal{D}}$ determinant with rows and columns indexed by $\mathcal{T}_{\mathcal{D}}^* = \{i_D, D \in \mathcal{D}, i_D \in \mathcal{I}_{\mathcal{D}}^*\}$ ordered in an arbitrary manner. For $i_D \in \mathcal{T}_{\mathcal{D}}^*$, the $i_D$-column is the vector of derivatives of $p(i(D))$, with respect to $\theta(j_C), j_C \in \mathcal{T}_{\mathcal{D}}^*$. From (2.19), straightforward differentiation shows that

$$(A.7) \qquad \frac{dp(i(D))}{d\theta(j_C)} = p(i(D)) \Big[ \delta_{(i_D)_C}(j_C) - \sum_{\substack{(l_H)_C = j_C \\ H \in \mathcal{E}_{\ominus}, l_H \in \mathcal{I}_E^*}} p(l(H)) \Big],$$

where

$$\delta_{(i_D)_C}(j_C) = \begin{cases} 1, & \text{if } (i_D)_C = j_C, \\ 0, & \text{otherwise.} \end{cases}$$

We note that the factor $p(i(D))$ is common to all components of the $i_D$ column and therefore the Jacobian is equal to

$$J = \det A \prod_{i_D \in \mathcal{T}_{\mathcal{D}}^*} p(i(D)),$$

where $A$ is the $d_{\mathcal{D}} \times d_{\mathcal{D}}$ matrix with entries

$$\delta_{(i_D)_C}(j_C) - \sum_{\substack{(l_H)_C = j_C \\ H \in \mathcal{E}_{\ominus}, l_H \in \mathcal{I}_E^*}} p(l(H)), \qquad i_D \in \mathcal{T}_{\mathcal{D}}^*, \ j_C \in \mathcal{T}_{\mathcal{D}}^*.$$

We also note that if $C$ is maximal, with respect to inclusion, for the row $r(j_C)$ of $A$ corresponding to $j_C \in \mathcal{T}_{\mathcal{D}}^*$, the only entry for which $\delta_{(i_D)_C}(j_C) \neq 0$ is the diagonal entry; therefore, for $C$ maximal, we can write

$$(A.8) \qquad r(j_C) = e(j_C) - \Big( \sum_{\substack{(l_H)_C = j_C \\ H \in \mathcal{E}_{\ominus}, l_H \in \mathcal{I}_E^*}} p(l(H)) \Big) \mathbf{1}^t,$$

where $e(j_C)$ is the $d_{\mathcal{D}}$-dimensional vector with components all equal to 0 except for the $j_C$ component, which is equal to 1, and $\mathbf{1}$ is the vector $d_{\mathcal{D}}$-dimensional vector with all its components equal to 1. If $C$ is not maximal, the entries corresponding to the columns $i(D)$ with $(i_D)_C = j_C$ (which implies that $C \subseteq D$) also have $\delta_{(i_D)_C}(j_C) = 1$. In order to eliminate, in the row $r(j_C)$ of $A$, the $\delta_{(i_D)_C}(j_C) = 1$ for $C$ strictly included in $D$, we replace $r(j_C)$ by

$$\tilde{r}(j_C) = r(j_C) + \sum_{\substack{(l_F)_C = j_C \\ F \supset C}} (-1)^{|F \setminus C|} r(l_F),$$

which yields, for any $j_C \in \mathcal{T}_{\mathcal{D}}^*$,

$$(A.9) \quad \tilde{r}(j_C) = e(j_C) - \left( \sum_{F \supseteq C, F \in \mathcal{D}} (-1)^{|F \setminus C|} \sum_{\substack{(l_H)_C = j_C \\ H \in \mathcal{E}_\ominus, l_H \in \mathcal{I}_E^*}} p(l(H)) \right) \mathbf{1}^t.$$

Clearly, for $C$ maximal, $\tilde{r}(j_C) = r(j_C)$. Moreover, the matrix $\tilde{A}$ obtained by replacing $\tilde{r}(j_C)$ by $r(j_C)$ has the same determinant, as $A$ and is equal to

$$\tilde{A} = I_{d_{\mathcal{D}}} - U \mathbf{1}^t,$$

where $U$ is the column vector

$$U = \left( \sum_{F \supseteq C, F \in \mathcal{D}} (-1)^{|F \setminus C|} \sum_{\substack{(l_H)_C = j_C \\ H \in \mathcal{E}_\ominus, l_H \in \mathcal{I}_H^*}} p(l(H)), j_C \in \mathcal{T}_{\mathcal{D}}^* \right).$$

It is well known that, for $\tilde{A}$ of that form, $\det \tilde{A} = 1 - \mathbf{1}^t U$; that is,

$$\det \tilde{A} = 1 - \sum_{j_C \in \mathcal{T}_{\mathcal{D}}^*} \left( \sum_{F \supseteq C, F \in \mathcal{D}} (-1)^{|F \setminus C|} \sum_{\substack{(l_H)_C = j_C \\ H \in \mathcal{E}_\ominus, l_H \in \mathcal{I}_H^*}} p(l(H)) \right)$$

$$= 1 - \sum_{H \in \mathcal{E}_\ominus, l_H \in \mathcal{I}_H^*} p(l(H)) \left( \sum_{F \subseteq H} \sum_{C \subseteq \ominus F} (-1)^{|F \setminus C|} \right)$$

$$= 1 - \sum_{H \in \mathcal{E}_\ominus, l_H \in \mathcal{I}_H^*} p(l(H)) \left[ \sum_{F \subseteq \mathcal{D} H} (-1)^{|F| - 1} \right],$$

where the last equality follows from the general fact that $\sum_{C \subseteq \ominus F} (-1)^{|F \setminus C|} = (-1)^{|F| - 1}$. Therefore, (4.14) for the general hierarchical model is now proved.

From (4.9) and Lemma 4.2, we see that in the particular case of graphical models, we have

$$\det \tilde{A} = p_\varnothing + \sum_{\substack{H \in \mathcal{U}_\ominus \\ i_H \in \mathcal{I}_H^*}} p(i(H)) + \sum_{\substack{H \notin \mathcal{U}_\ominus \\ i_H \in \mathcal{I}_H^*}} p(i(H))$$

$$- \sum_{\substack{H \in \mathcal{E}_{\ominus} \\ l_H \in \mathcal{I}_H^*}} p(l(H)) \left( \sum_{F \subseteq_{\mathcal{D}} H} (-1)^{|F|-1} \right)$$

$$= p_{\varnothing} + \sum_{\substack{H \in \mathcal{U}_{\ominus} \\ i_H \in \mathcal{I}_H^*}} \left[ p(i(H)) \left( 1 - \sum_{F \subseteq_{\mathcal{D}} H} (-1)^{|F|-1} \right) \right]$$

$$+ \sum_{\substack{H \notin \mathcal{U}_{\ominus} \\ i_H \in \mathcal{I}_H^*}} \left[ p(i(H)) \left( 1 - \sum_{F \subseteq_{\mathcal{D}} H} (-1)^{|F|-1} \right) \right]$$

$$= p_{\varnothing} + \sum_{\substack{H \in \mathcal{U}_{\ominus} \\ i_H \in \mathcal{I}_H^*}} \left[ p(i(H)) \left( 1 - \sum_{F \subseteq_{\mathcal{D}} H} (-1)^{|F|-1} \right) \right]$$

$$= p_{\varnothing} - \sum_{\substack{H \in \mathcal{U}_{\ominus} \\ i_H \in \mathcal{I}_H^*}} [p(i(H)) a(H)],$$

which proves (4.15).

H. Massam
J. Liu
Department of Mathematics
    and Statistics
York University
Toronto, Ontario M3J 1P3
Canada
E-mail: massamh@yorku.ca
        jnliu@yorku.ca
URL: http://www.math.yorku.ca/massamh

A. Dobra
Department of Statistics
University of Washington
Seattle, Washington 98195
USA
E-mail: adobra@u.washington.edu
URL: http://www.stat.washington.edu/adobra